\documentclass[10pt,a4 paper]{article}
\usepackage{amssymb}
\usepackage{amsmath}
\usepackage{amsthm}
\usepackage{amsfonts}
\usepackage{amstext}
\usepackage{cite}
\usepackage{color}
\usepackage{hyperref} 

\newtheorem{theorem}{Theorem}[section]
\newtheorem{defff}[theorem]{Definition}

\newtheorem{rem}[theorem]{Remark}
\newtheorem{lee}[theorem]{Lemma}
\newtheorem{col}[theorem]{Corollary}

\numberwithin{equation}{section}

\newcommand{\ds}{\displaystyle}

\begin{document}
	
\title{Parameter analysis in continuous data assimilation for three-dimensional  Brinkman-Forchheimer-extended Darcy model}
\author{  D\'ebora A.F. Albanez$^{\lowercase{a}}$,  Maicon José Benvenutti$^{\lowercase{b}}$, }
\maketitle

\begin{center}

$^a$ Departamento Acadêmico  de Matem{\'a}tica \; -- \;
Universidade Tecnol{\'o}gica Federal do Paran{\'a}, \\
86300-000 Corn{\'e}lio Proc{\'o}pio, PR -- Brasil, \\corresponding author: deboraalbanez@utfpr.edu.br\\ 
$^b$ Departamento de Matemática \; -- \;
Universidade Federal de Santa Catarina, \\
89036-004 Blumenau, SC -- Brasil,\\ m.benvenutti@ufsc.br\\

\end{center}

\begin{abstract}
	In this paper, we study analytically the long-time behavior of three-dimensional  Brinkman-Forchheimer-extended Darcy model, in the context that the parameters related to the damping nonlinear term are unknown. This work is inspired by the approach firstly introduced for two-dimensional Navier-Stokes equations by Carlson, Hudson and Larios. We show estimates in $L^2$ and $H^1$ for large-time error between the true solution and the assimilated solution, which is constructed with the unknown damping parameters and observational measurements obtained continuously in time from a continuous data assimilation technique proposed by Azouani, Olson and Titi.
	
	
	\noindent\textbf{AMS 2020 MSC:} 93C20, 76D05, 37C55, 35Q35, 35Q30, 	34D06.
	\vspace{0.2cm}
	
	
	\vspace{0.2cm}
\end{abstract}


\section{Introduction}


In an attempt to determinate a mathematical model as accurate as possible which represents the dynamics of any physical phenomena, parameters related to natural features of the involved elements of the system  frequently arise, for instance, the external force and the kinematic viscosity of an incompressible fluid on analysis of its velocity through Navier-Stokes equations, or thermal viscosity if the temperature of the fluid is an object of interest. Thus, the development of data-driven techniques for the purpose of estabilishing accurate values of the parameters, namely parameter learning, becomes a fundamental point of the analysis of the dynamical system, in order to improve the model.

Recently, parameter estimation algorithms were applied to two dimensional Navier-Stokes equations (see \cite{Carlson1}) and three dimensional Lorenz system (see \cite{Carlson2}). In a sensitivity-type analysis, results in \cite{Carlson1} and \cite{Carlson2} exhibited the large-time error between the true solution of the model and the assimilated solution due to the deviation between the approximate and physical parameters.

Based on ideas of \cite{Carlson1} and \cite{Carlson2}, we consider the three-dimensional Brinkman-
Forchheimer-extended Darcy model, also named as three-dimensional Navier-Stokes with damping equations 
\begin{equation}
	\left\{
	\begin{array}{l}
		\displaystyle \frac{\partial u}{\partial t}+(u\cdot\nabla)  u- \nu \Delta u +\nabla p +a|u|^{2\alpha}u
		=f,\\ \\
		\nabla\cdot u = 0,
	\end{array}
	\right.  \label{1aux}
\end{equation}%
where  $u=(u_{1}(x,t),u_{2}(x,t),u_{3}(x,t))$ is the spatial  velocity field,  $p=p(x,t)$ is a  scalar pressure field, $f=f(x,t)$ is a given external force and $\nu >0$ is the kinematic viscosity. The parameters of interest in this work are $\alpha>1$ and $a>0$, i.e., the coefficients related to the damping nonlinear term $a|u|^{2\alpha}u$, which we suppose being both unknown. Notice that in the limit case $a=0$, we obtain the classical Navier–Stokes system. Equation (\ref{1aux}), as well as  most fluid models in porous media recently studied, arised from Darcy's Law, that describes a proportional relation between the instantaneous rate of discharge through a porous medium,  the viscosity of the fluid, and the pressure drop over a given distance.

Darcy’s law is generally valid for flows with Reynold's number $Re\leq 1$, i.e., laminar flows. For the treatment of  Darcy´s Law deviation cases, a more suitable model is obtained by coupling a quadractic term in seepage velocity to account for the increased pressure drop, namely the Darcy-Forchheimer equation (see \cite{Amao}):
$$\frac{\mu}{k}v+\beta\rho|v|^2v=-\nabla p, $$
where $\mu$ is the kinematic viscosity, $k$ is the permeability of the porous medium, $\beta$ is the inertial factor and $\rho$ is the density of the fluid flowing through the medium. For the system \eqref{1aux}, $a|u|^{2\alpha}u$ is a drag term related to the pore dimension, shape and porosity (see \cite{{Giorgi}}, \cite{HSU}, \cite{Nield}, \cite{Nield2} and \cite{{Whitaker}}). Although the most widely values used for  $\alpha$ are $0$, $\frac{1}{2}$ and $1$, extrapolations to others real numbers have appeared in literature and the suitable range for modeling purpose is still source of uncertainty (see \cite{Amao}, \cite{Markowich}, \cite{Mei}, \cite{Mei2} and \cite{Vafai}). Due to the mathematical restriction related with existence and uniqueness of solutions, in this paper we will consider the case $\alpha >1$.

 The aim of this work is to present a parameter analysis with respect to errors on the damping term parameters $\alpha$ and $a$, based on the technique of a continuous data assimilation algorithm (named here as AOT algorithm) proposed in \cite{Azouani}. This algorithm, initially applied to two-dimensional Navier-Stokes equations, was designed to work for general linear and nonlinear dissipative dynamical systems, based on a feedback control that works inserting the large scale observations into the physical model through of a linear interpolant operator constructed from these observational measurements. This feedback control  is basically used for relaxing the solution of the constructed system towards the real-state solution of the original model. In past recently years, it has been analyzed for several important 1D, 2D and 3D physical models (see 
 \cite{Azouani2}, \cite{Albanez1}, \cite{Albanez2}, \cite{Biswas},\cite{Biswas2}, \cite{Farhat1}, \cite{Farhat2}, \cite{Farhat4}, \cite{Jolly1}, \cite{Jolly2}, \cite{Jolly3}).

  We prove that under certain hypothesis, one can estimate $u$ asymptotically in $H^{1}$-norm for observable interpolants accurate enough in $L^{2}$-norm (see Theorems \ref{CDA3a1.b.c.d} and \ref{CDA3a1.b.c.d.e}). For this purpose, we consider the assimilated system given by
\begin{equation}
	\left\{
	\begin{array}{l}
		\displaystyle \frac{\partial w}{\partial t}+(w\cdot\nabla)  w- \nu \Delta w +\nabla p +b|w|^{2\beta}w
		=f+ \eta(I_{h}(u)-I_{h}(w)),\\ \\
		\nabla\cdot w = 0,
	\end{array}
	\right.  \label{1qw.n} 
\end{equation}%
where $\eta>0$ is the relaxation (nudging) parameter, $u$ is the solution of  (\ref{1aux}), $b>0$ is a guess for $a$ and $\beta>1$ is a guess for $\alpha$, $I_{h}$ is a linear interpolation operator satisfying certain conditions (see (\ref{0in1}) below) and $h>0$ is a parameter related with spatial resolution of the operator.



We consider the systems (\ref{1aux}) and (\ref{1qw.n}) subjected to the periodic boundary condition  on a box domain $\Omega=[0,l]^{3}$, namely, 
\begin{equation}
	u(x,t)=u(x+le_{i},t),\,\, \forall\, (x,t)\in \mathbb{R}^{3}\times \lbrack
	0,T ), \label{000}
\end{equation} 
where $e_{1}$,\,$e_{2}$ and $e_{3}$ are the canonical basis of $\mathbb{R}^{3}$ and $l>0$ the fixed period. For the assimilated operator, we assume that
\begin{equation} \label{0in00}
	I_{h}:(H^{j}([0,l]^{3})^{3} \longrightarrow (L^{2}([0,l]^{3})^{3}
\end{equation}
is a linear operator where there exist dimensionless constants $c_0$ and $c_1$ such that $c_1=0$ and $c_0>0$ if $j=1$, $c_0>0$ and $c_1>0$ if $j=2$, and the following inequality is satisfied: 
\begin{equation}\label{0in1} 
	\|I_{h}(g)-g\|^{2}_{L^{2}}\leq c_{o}h^{2}\|\nabla g\|^{2}_{L^{2}}+ c_{1}h^{4}\|\Delta g\|^{2}_{L^{2}}, \,\, \forall \,   g \in (H^{j}([0,l]^{3})^{3},
\end{equation}



The outline of this paper is as follows.: In Section \ref{sb1}, we start stating the mathematical
	setting, notations and classical inequalities used for obtaining all the results. Furthermore, weak and strong solutions for systems of interest  \eqref{1aux} and \eqref{1qw.n} are defined. The known results of global well-posedness for these systems are also presented.
	
In Section \ref{yhdbh}, we enunciate the results obtained: Theorems \ref{CDA3a1.b}, \ref{CDA3a1.b.c.d} and \ref{CDA3a1.b.c.d.e} assert that considering  the observational measurements of the system \eqref{1aux} jointed with the technique of AOT algorithm,  the true solution of the system can be recovered with only approximated values of the parameters related to damping term, namely $a$ and $\alpha$,  except for a remaining error  of the difference between the real parameter and approximated inserted one.

In Section \ref{sec4}, estimates in $L^2$ and other apropriate spaces for the solution of the system \eqref{1aux} are performed, as well as for space and time derivatives; Such results are applied on proofs of the theorems established in Section 3.  Finally, Sections \ref{sec5}, \ref{sec6} and \ref{sec7} contain the proofs of Theorems \ref{CDA3a1.b}, \ref{CDA3a1.b.c.d} and \ref{CDA3a1.b.c.d.e}, respectively.


\section{Funtional setting and results}
\label{sb1}


\subsection{Functions spaces, functionals and inequalities}


For $\Omega=[0,l]^3$ and $m \in \mathbb{N}$, let  $L^{m}(\Omega)$ be the usual $l$-periodic Lebesgue $m$-integrable space,  $H^{m}(\Omega)$ the usual $l$-periodic trigonometric Sobolev space (see \cite{Rosa}) and $H$ and $V$ the vector subspaces given by
\begin{align}
	H= \left\{u \in (L^2(\Omega))^{3}\,;\,\, \nabla \cdot u= 0\right\}, \,\,\,
	V= \left\{u \in (H^{1}(\Omega))^{3}\,;\,\, \nabla \cdot u= 0\right\}\nonumber
\end{align}
and endowed with the product structures. 
For each fixed $\alpha  > 1$, let us also define the Banach space $Y_{\alpha}=V\cap (L^{2\alpha+2}(\Omega))^{3}$ and, then $Y_{\alpha}'=V'+(L^{\frac{2\alpha+2}{2\alpha+1}}(\Omega))^{3}$. 


Let $A: Y_{\alpha}\rightarrow Y_{\alpha}'$ be the operator defined by $$\langle A(u),v\rangle_{Y_{\alpha}',Y_{\alpha}}=\int_{\Omega} \nabla u\cdot \nabla v\, dx,$$ and, for each $0 \leq \gamma \leq \alpha $,  $G_{\gamma}:  Y_{\alpha} \rightarrow Y_{\alpha}'$ the operator given by $$\langle G_{\gamma}(u),v\rangle_{Y_{\alpha}', Y_{\alpha}} =\int_{\Omega} |u|^{2\gamma}u\cdot v\, dx.$$ We also define  the bilinear operator $B:Y_{\alpha}\times Y_{\alpha}\rightarrow Y_{\alpha}'$ as 
$$\langle B(u,v),w \rangle_{Y_{\alpha}',Y_{\alpha}}= \int_{\Omega}(u\cdot \nabla v) \cdot w\, dx.$$
For $u,v, w \in Y$, the term $B$ has the property
\begin{equation}
	\langle B(u,v),w \rangle_{Y_{\alpha}',Y_{\alpha}}=-\langle B(u,w),v\rangle_{Y_{\alpha}',Y_{\alpha}}, \label{fdr1}
\end{equation}
and hence 
\begin{equation}
	\langle B(u,w),w\rangle_{Y_{\alpha}',Y_{\alpha}}=0 \label{fdr2}.
\end{equation}
Moreover, we have that
\begin{equation}
	-\langle B(u,u)-B(v,v),u-v\rangle_{Y_{\alpha}',Y_{\alpha}}=\frac{1}{2} \langle B(u-v,u-v),u+v \rangle_{Y_{\alpha}',Y_{\alpha}}. \label{fdr3}
\end{equation}


Let $\mathcal{P}:(L^2(\Omega))^{3}\rightarrow H$ be the classical Helmholtz-Leray orthogonal projection (see \cite{Rosa}). If $u,v, \in V\cap (H^{2}(\Omega))^{3}$, then 
\begin{equation}
	Au=-\mathcal{P}(\Delta u)= -\Delta u \,\,\mbox{ and }\,\, B(u,v)= \mathcal{P}(u\cdot \nabla v), \nonumber
\end{equation}
and an equivalent norm in $V\cap (H^{2}(\Omega))^{3}$ is given by
\begin{eqnarray}
	\left(\|u\|^{2}_{L^{2}}+ \|Au\|^{2}_{L^{2}}\right)^{\frac{1}{2}}, \nonumber
\end{eqnarray}
where, hereafter, the domain $\Omega=[0,l]^3$ is omitted in the expressions involving norms.


We recall some particular three-dimensional cases of the Gagliardo-Nirenberg inequality  (see \cite{Friedman} and \cite{Temam2}):
\begin{equation}
	\left\{
	\begin{array}{ll}
		\displaystyle \|  g\|_{L^{\infty}} \leq C_{\infty}\left(\| \nabla  g\|_{L^{2}}^{\frac{1}{2}}\| Ag\|_{L^{2}}^{\frac{1}{2}}+\frac{1}{l^{\frac{3}{2}}}\|g\|_{L^{2}}\right), \,\, \forall\, g \in V\cap (H^{2}(\Omega))^{3}; & \\&\\ \displaystyle
		\| g\|_{L^{p}}\leq C_{p}\left(\|\nabla g\|^{\frac{6+p}{2p}}_{L^{2}}\|A g\|^{\frac{p-6}{2p}}_{L^{2}} +\frac{1}{l^{\frac{3p-6}{2p}}}\|g\|_{L^{2}}\right), \,\, \forall\, 6 < p<\infty \,\, \mbox{ and } \,\, g \in V\cap (H^{2}(\Omega))^{3}; &  \\&\\ \displaystyle
		\| g\|_{L^{p}}\leq C_{p}\left(\|g\|^{\frac{6-p}{2p}}_{L^{2}}\|\nabla g\|^{\frac{3p-6}{2p}}_{L^{2}} +\frac{1}{l^{\frac{3p-6}{2p}}}\|g\|_{L^{2}}\right), \,\, \forall\, 2<p\leq 6 \,\, \mbox{ and } \,\, g \in V, &  
	\end{array}%
	\right.  \label{gn}
\end{equation}%
where  $C_{p}$ are dimensionless constants. 


Furthermore, we use the following inequality to deal with the nonlinear damping term (see \cite{Li}):
\begin{equation}
	(|x|^{\gamma}x - |y|^{\gamma}y) \cdot (x-y) \geq \frac{1}{2}|x-y|^{2}(|x|^{\gamma}+|y|^{\gamma})\:\:\,\, \forall\, x,\, y \in \mathbb{R}^{3} \,\, \mbox{ and }\,\, \gamma \geq 0, \label{dam01}
\end{equation}
and, by Mean Value Theorem, there exists a dimensionless constant $\kappa(\gamma)>0$ such that
\begin{equation}
	||x|^{\gamma}x-|y|^{\gamma}y|\leq \kappa(\gamma)|x-y|(|x|+|y|)^{\gamma}\:\:\,\, \forall\, x,\, y \in \mathbb{R}^{3} \,\, \mbox{ and }\,\, \gamma \geq 0. \label{dam02}
\end{equation} 


\subsection{Weak and strong solutions}


In order to define and analyze weak and strong solutios for equation (\ref{1aux}), we rewrite it using functional settings as 
\begin{equation}
	\left\{
	\begin{array}{l}
		\displaystyle \frac{du}{dt}+\nu Au+B(u,u)+aG_{\alpha}(u) =\mathcal{P}(f), \\ \\
		u(0)=u_{0},
	\end{array}
	\right. \label{eq111}
\end{equation}%
and (\ref{1qw.n}) as
\begin{equation}
	\left\{
	\begin{array}{l}
		\displaystyle \frac{dw}{dt}+\nu Aw+B(w,w)+bG_{\beta}(w)=\mathcal{P}(f)+\eta\mathcal{P} (I_{h}(u)-I_{h}(w))\\ \\
		w(0)=w_{0}.
	\end{array}
	\right. \label{eq1112}
\end{equation}


\begin{defff}[Weak solution]
	Suppose $\alpha > 1$ and $a>0$. Let $f\in L^{2}((0,T),H)$ and $u(0)=u_{0}\in H$. A local weak solution for  system (\ref{eq111}) is a function $u \in L^{\infty}((0,T),H)\cap L^{2}((0,T),V)\cap L^{2\alpha+2}((0,T),(L^{2\alpha+2}(\Omega))^{3})$ such that satisfies (\ref{eq111})
	in $L^{1}((0,T), Y_{\alpha}')$. We say that $u$ is a global weak solution if $u$ is local weak solution for each $T>0$.
\end{defff}
\begin{rem} \label{remouio}
	If $u$ is a weak solution for  system (\ref{eq111}), then 
	\begin{equation}
		\frac{du}{dt} \in L^{2}((0,T),V')+L^{\frac{2\alpha+2}{2\alpha+1}}\left((0,T),(L^{\frac{2\alpha+2}{2\alpha+1}}(\Omega))^{3}\right), \nonumber
	\end{equation}
	$u \in  C([0,T),H)$  and  $\displaystyle \left\langle\frac{du}{dt},u \right\rangle_{Y_{\alpha}',Y_{\alpha}}=\frac{d}{dt}\|u\|^{2}_{L^{2}}$ (see \cite{Li2} and \cite{Pardo}).
\end{rem}


\begin{defff}[Strong solution]
	Suppose $\alpha >1$ and $a>0$. Let $f\in L^{2}((0,T),H)$ and $u(0)=u_{0}\in V$. A local strong solution for  system (\ref{eq111}) is a weak  solution such that $u \in L^{\infty}((0,T),V)\cap L^{2}((0,T),(H^{2}(\Omega))^{3})$.
	We say that $u$ is a global strong solution if $u$ is local strong solution for each $T>0$.
\end{defff}
\begin{rem}
	If $u$ is a strong solution for  system (\ref{eq111}), then 
	\begin{equation}
		\frac{du}{dt} \in L^{2}\left((0,T),H\right) + L^{\frac{2\alpha+2}{2\alpha+1}}\left((0,T),(L^{\frac{2\alpha+2}{2\alpha+1}}(\Omega))^{3}\right)  \nonumber
	\end{equation}
	and   $u \in  C([0,T),V)$  (see \cite{Temam}).
\end{rem}


\begin{theorem}[Global existence and uniqueness of weak and strong solutions] \label{ExiUniBard} Suppose $\alpha > 1$, $a>0$ and $f\in L^{2}_{loc}(\mathbb{R}_{+},H)$. If $u(0)=u_{0}\in H$, then
	the system (\ref{eq111}) has an unique global weak solution, which is continuously dependent on the initial data in the $H$-norm. Furthermore, if $u(0)=u_{0}\in V$, then the global weak solution is strong and it is also continuously dependent on the initial data in the $V$-norm.
\end{theorem}
The proofs of the above result can be found \cite{Markowich} (see also
\cite{Cai},  \cite{Li}, \cite{Li2}, \cite{Zhou}, \cite{Lu}, \cite{Zhong}  and \cite{Zhong2}). \\

Concerning to the assimilated system (\ref{eq1112}), similar definitions of weak and strong solutions can be stated. Likewise, there are the following results:
\begin{theorem}[Existence and uniqueness of weak and strong solutions for (\ref{eq1112}) with $c_1=0$] \label{ExiUniBard3a1}
	Suppose $\alpha, \,\, \beta > 1$, $a, \,\,b>0$ $f\in L^{2}_{loc}(\mathbb{R}_{+},H)$, $u$ a global weak solution of (\ref{eq111}) and $I_{h}$ a linear operator that satisfies (\ref{0in1}) with $c_{1}=0$. If  $w(0)=w_{0}\in H$, then the system  (\ref{eq1112}) has an unique global weak solution, which is continuously dependent on the initial data in the $H$-norm. Furthermore, if  $w(0)=w_{0}\in V$, and $u$ is a global strong solution of (\ref{eq111}), then the global weak solution of (\ref{eq1112}) is strong and also continuously dependent on the initial data in the $V$-norm.
\end{theorem}

\begin{theorem}[Existence and uniqueness of strong solutions for (\ref{eq1112}) with $c_{1} \geq 0$] \label{ExiUniBard3a2}
	Suppose $\alpha, \,\, \beta > 1$, $a, \,\,b>0$ $f\in L^{2}_{loc}(\mathbb{R}_{+},H)$, $u$ a global strong solution of (\ref{eq111}) and $I_{h}$ a linear operator that satisfies (\ref{0in1}). If  $w(0)=w_{0}\in V$, and $\nu^{2}> \eta^{2}c_{1}h^{4}$, then the system  (\ref{eq1112}) has an unique global strong solution and there is continuous dependence on the initial data in the $V$-norm. 
\end{theorem}


The proof of above results can be also found in \cite{Markowich}. 


\section{Statements of main theorems}
\label{yhdbh}


In this work, we are considering dimensional equations in (\ref{eq111}) and (\ref{eq1112}) and all the results and inequalities presented here are in a correct balance of units. Since $\alpha$ and $\beta$ are dimensionless constants, we have that the unit of 
$a$ is  $\displaystyle \frac{\left[\mbox{time}\right]^{2\alpha-1}}{\left[\mbox{length}\right]^{2\alpha}}$, while $b$ is $\displaystyle \frac{\left[\mbox{time}\right]^{2\beta-1}}{\left[\mbox{length}\right]^{2\beta}}$. Thus, for suitable balance of the equations \eqref{eq111} and \eqref{eq1112}, we write
\begin{align}
	a=\tilde{a} \left(\frac{l^{2\alpha-2}}{\nu^{2\alpha-1}}\right) \,\,\,\, \mbox{ and } \,\,\,\, b=\tilde{b} \left(\frac{l^{2\beta-2}}{\nu^{2\beta-1}}\right), \label{4575409}
\end{align}
where $\tilde{a}$ and $\tilde{b}$ are dimensionless constants. So our estimates will be given in term of differences involving $\tilde{a}$ and $\tilde{b}$. We observe that (\ref{4575409})  is in agreement with expressions commonly used in the literature for the most usual cases $\alpha$ equals to $0$, $\frac{1}{2}$ and $1$
(see, for example, \cite{Ingham}, \cite{Joseph} and \cite{Skjetne}).


We start with a result in the context of $L^{2}$-norm. In order to express the estimates properly, we consider  a upper bound in $V$-norm for initial data of the physical system (\ref{eq111}) given in (\ref{estM}). Besides, we have the restrictions $ 1<\alpha,\,\beta<3$ and $c_{1}=0$ in (\ref{0in1}).


\begin{theorem} \label{CDA3a1.b} Let $f\in L^{\infty}(\mathbb{R}_{+};H)$, $1<\alpha, \,\beta<3$ and $a, \, b>0$ given. Suppose that the linear interpolation operator $I_h$ satisfies (\ref{0in1}) with $c_0>0$ and $c_1=0$. Consider $\eta$ and $h$ large and small enough, respectively, such that
	\begin{align}
		\eta&> \frac{8(\beta-1)}{\beta b^{\frac{1}{\beta-1}}\nu^{\frac{\beta}{\beta-1}}} \,\,\mbox{ and }\,\, \nu > 4\eta c_{o}h^{2}, \label{12ee2} 
	\end{align}
	Moreover, let $u$ be a global strong solution of (\ref{eq111}), $M >0$ such that
	\begin{align}
		\|\nabla u(0)\|_{L^{2}}^{2}+ \frac{1}{l^{2}} \|u(0)\|_{L^{2}}^{2}\leq M. \label{estM}
	\end{align}
	and $w$ a global weak solution of  (\ref{eq1112}). If $1 < \alpha,\, \beta<2$, we have
	\begin{align}
		\|w(t)-u(t)\|_{L^{2}}^{2} \leq e^{-\frac{\eta}{8}t}\|w(0)-u(0)\|_{L^{2}}^{2} &+|\alpha-\beta|^{2}\left[32\tilde{a}^{2}\frac{\nu^{2}}{\eta l^{4}}M_{1} + \frac{64\tilde{a}^{2}C^{12}_{6}}{(2-\max\{\alpha,\,\beta\})^{2}}A_{0}\right]
		\nonumber\\&+ \left|\tilde{a}-\tilde{b}\right|^{2}\left[2\frac{\nu^{2}}{\eta l^{4}}M_{1}+ 2C^{10}_{6}A_{0} \right],
	\end{align}%
for all $t\geq 0$, where
	\begin{align}
		A_{0}&= \frac{l^{2}(\eta l^{2}+2\nu)}{\eta^{2}\nu^{7}}\left(M_{2}+\frac{1}{l^{2}}M_{1}\right)^{5}
	\end{align}
	and $M_{1}$ and $M_{2}$ are constants given in Corollary \ref{da0101.4565.1} and uniform estimates for norms of $u$ and $\nabla u$ in $L^{2}$, respectively.
	
Furthermore, if $2 \leq \alpha <3$ or $2 \leq \beta <3$, we have
	\begin{align}
		\|w(t)-u(t)\|_{L^{2}}^{2} \leq e^{-\frac{\eta}{8}t}\|w(0)-u(0)\|_{L^{2}}^{2} &+|\alpha-\beta|^{2}\left[2^{9}\tilde{a}^{2} \frac{\nu^{2} }{\eta^{2} l^{4}}M_{1}+
		\frac{2^{22}\tilde{a}^{2}C^{2}_{6} C^{14}_{\frac{42}{5}}}{(3-\max\{\alpha,\,\beta\})^{2}}A_{1}\right]
		\nonumber\\&+ \left|\tilde{a}-\tilde{b}\right|^{2}\left[2^{5}\frac{\nu^{2}}{\eta^{2}l^{4}}M_{1}+2^{16}C^{14}_{\frac{42}{5}}A_{1}\right],
	\end{align}%
for all $t\geq 0$, with
	\begin{align}
		A_{1}&= \frac{l^{8}}{\nu^{10}}\left(\frac{1}{\nu}+\frac{2}{\eta l^{2}}\right)\left[\frac{1}{\nu}M^{7}_{2}+ \frac{4}{\eta \nu^{\frac{3\alpha-2}{\alpha-1}} a^{\frac{1}{\alpha-1}}} M^{7}_{2} +\frac{8}{\eta \nu^{2}}\|f\|^{2}_{L^{\infty}_{t}L^{2}}M^{6}_{2}+\frac{2}{\eta l^{16}}M_{1}^{7}\right]. 
	\end{align}
\end{theorem}


\begin{rem}
	Note that the relaxation parameter $\eta$ does not appear on $A_{1}$ term $\frac{l^{8}M_{2}^{7}}{\nu^{12}}$, while it is present in fraction denominator of each term of $A_{0}$. Therefore, for the case $1 < \alpha, \,\beta<2$, it is possible that the error of approximation $\|w(t)-u(t)\|_{L^{2}}$ be small enough in the asymptotic sense by choosing $\eta$ large enough (and consequently $h$ small enough, namely, $I_{h}$ accurate enough), regardless of the parameter errors $|\alpha-\beta|$ and $\left|\tilde{a}-\tilde{b}\right|$. 
\end{rem}

\begin{rem}
	Since $u$ is a strong solution to (\ref{eq111}), by using (\ref{0in1}) with $c_{1}=0$ and (\ref{12ee2}), we obtain  $$\|I_{h}(u(t))-u(t)\|^{2}_{L^{2}}\leq \frac{\nu}{4\eta}M_{2},\,\,\text{for all}\:\:  t \geq 0,$$ and thus $I_{h}(u)$ can also be used directly as an approximation to $u$. For many cases, combining approaches to obtain a better approximation is more suitable. For instance, if $I_{h}$ is the projection onto low Fourier modes, we can consider  $$u(t) \approx I_{h}(u(t))+(I-I_{h})(w(t)),$$ with $I$ the identity operator. Here, the low modes values are extracted from $I_{h}(u)$  while the high modes values are from $w$.
\end{rem}



If we consider only strong solutions of  (\ref{eq1112}), we have results in $V$-norm. In this case, we need also an estimate for the initial data of the assimilated system (\ref{eq1112}), as given in (\ref{dgdbrt}). Besides, we also have the restriction $ 1<\alpha, \,\beta<2$, but $c_{1} \geq 0$.

\begin{theorem}\label{CDA3a1.b.c.d} Suppose $f\in L^{\infty}(\mathbb{R}_{+};H)$, $1<\alpha, \beta <2$ and $a, \, b>0$. Let $u$ be a global strong solution of (\ref{eq111}). Consider the linear interpolation operator $I_h$ satisfiying (\ref{0in1}) and $M >0$ satisfying (\ref{estM}), as well as $M_{3}$ and $\tilde{Z}_{1}$  constants given in (\ref{dfbdfr54}) and Corollary \ref{da0101.4565.1}.	Consider also $\eta$ and $h$ large and small enough, respectively, such that  (\ref{12ee2}) is satisfied and
	\begin{align}
		\eta>\frac{32\eta^{2}c_{0}h^{2}}{\nu}+4\tilde{Z}_{1} + 2^{14}\kappa^{2}(2\beta)C^{2}_{6}C^{4\beta}_{6\beta}\frac{l}{\nu}M_{3}  \,\,\mbox{ and }\,\, \nu^{2}> \frac{32c_{1}}{7}h^{4}\eta\frac{(\nu+8\eta l^{2})}{l^{2}}. \label{45drht}
	\end{align}
Moreover, let $w$ be a global strong solution of  (\ref{eq1112}) with
	\begin{align}
		\|\nabla w(0)\|_{L^{2}}^{2}+ \frac{1}{l^{2}} \|w(0)\|_{L^{2}}^{2}\leq M. \label{dgdbrt}
	\end{align}
	Then, for $B$, $C$ and $D$ constants given in (\ref{zzzer1}), (\ref{zzzer2}) and (\ref{zzzer3}), we have for all $t\geq 0$,
\begin{eqnarray}
	\|\nabla (w(t)-u(t))\|_{L^{2}}^{2}+ \frac{1}{l^{2}}\|w(t)-u(t)\|_{L^{2}}^{2} & \leq & Be^{-\frac{\eta}{8}t}\left(\|\nabla (w(0)-u(0))\|_{L^{2}}^{2}+ \frac{1}{l^{2}}\|w(0)-u(0)\|_{L^{2}}^{2}\right)\nonumber \\
	& + & C|\alpha-\beta|^{2} +D|\tilde{a}-\tilde{b}|^{2}. \label{67867ygb}
\end{eqnarray}

\end{theorem}


\begin{rem}
	Observe that $\tilde{Z}_{1}$ given in (\ref{dfbdfr54}) depends on constant $H$ appeared in (\ref{uioudrt5e5}), which in turn contains $\eta$ only on fraction denominators, that implies the same for 
$\tilde{Z}_{1}$. However, inequality (\ref{45drht}) can always be obtained for $\eta$ and $h$ large and small enough, respectively.
\end{rem}


Once again, note that in $C$ and $D$ given (\ref{zzzer2}) and (\ref{zzzer3}), the term  $\frac{l^{4}M_{2}^{5}}{\nu^{8}}$ has no relaxation term $\eta$ on its denominator. In the next result, we present an estimate where $\eta$ appears on fraction denominator of each term that multiplies $|\alpha-\beta|^{2}$ and $|\tilde{a}-\tilde{b}|^{2}$, thus, theoretically, being possible to make the error of approximation $\|w(t)-u(t)\|_{H^{1}}$ be small enough, regardless of the parameter errors, by having $\eta$ large enough. We restrict the analysis to the case where $f\in L^{\infty}(\mathbb{R}_{+};H)$ with $f_{t}\in L^{\infty}(\mathbb{R}_{+};H)$. Furthermore, the estimates are given from the time $\frac{2l^{2}}{\nu}$.


\begin{theorem}\label{CDA3a1.b.c.d.e} Let $f\in L^{\infty}(\mathbb{R}_{+};H)$ with $f_{t}\in L^{\infty}(\mathbb{R}_{+};H)$,  $u$  the global strong solution of (\ref{eq111}), with $1<\alpha<2$ and $a>0$. Consider $I_{h}$ a linear operator that satisfies (\ref{0in1})  and $w$ a global strong solution of (\ref{eq1112}) with  $1<\beta <2,\,b>0$ that satisfies (\ref{dgdbrt}). Also, suppose $M >0$ satisfying (\ref{estM}), $\eta$ large enough and $h$ small enough such that (\ref{12ee2}) and (\ref{45drht}) are valid. Then, for $B$, $\tilde{C}$ and $\tilde{D}$ constants given in (\ref{zzzer1}), (\ref{zzzer2.2}) and (\ref{zzzer3.3}), we have
	\begin{eqnarray}
	\:\:\:\:\:\|\nabla (w(t)-u(t))\|_{L^{2}}^{2}&\!\!\!\!\!\! + &\!\!\! \frac{1}{l^{2}}\|w(t)-u(t)\|_{L^{2}}^{2} \nonumber\\
		 &\leq &  Be^{-\frac{\eta}{8}\big(t- \frac{2l^{2}}{\nu}\big)}\left(\left\|\nabla (w\left(\tfrac{2l^{2}}{\nu}\right)-u\left(\tfrac{2l^{2}}{\nu}\right)\right\|_{L^{2}}^{2}+ \tfrac{1}{l^{2}}\left\|w\left(\tfrac{2l^{2}}{\nu}\right)-u\left(\tfrac{2l^{2}}{\nu}\right)\right\|_{L^{2}}^{2}\right)\nonumber\\
		 & + &  \tilde{C}|\alpha-\beta|^{2} +\tilde{D}|\tilde{a}-\tilde{b}|^{2},  \label{67867ygb.e}
	\end{eqnarray}
for all $t \geq \ds\frac{2l^{2}}{\nu}.$
\end{theorem}


\section{Estimates to the system (\ref{eq111})}\label{sec4}


Henceforth we present auxiliary results that will be useful in the proofs of Theorems \ref{CDA3a1.b} and \ref{CDA3a1.b.c.d} and they are based on energy-type estimates. To overcome eventual lack of regularities, these estimates are initially obtained for approximate solutions  coming from the Galerkin's procedure. Then, via a limit process, they are also obtained for the exact solutions. Since this procedure is standard, we will present the estimates directly on the exact system. 


\begin{lee} \label{da0101}
	Suppose $f\in L^{\infty}(\mathbb{R}_{+};H)$, $\alpha>1$, $a>0$ and 
	let $u$ be a global strong solution of (\ref{eq111}). Consider $K$  given by
	\begin{align}
		K&=\frac{l^{2}}{\nu}\|f\|^{2}_{L^{\infty}_{t}L^{2}}+ \frac{4\nu^{\frac{\alpha+1}{\alpha}}}{a^{\frac{1}{\alpha}}l^{\frac{2-\alpha}{\alpha}}}. \label{ct01} 
	\end{align}
	Then, we have the following estimates:
	
		\begin{equation}\label{jhgtn01}
		\|u(t)\|_{L^{2}}^{2} \leq e^{\frac{-2\nu t}{l^{2}}}\|u(0)\|^{2}_{L^{2}}+  \frac{l^{2}}{\nu}K,\,\, \forall\, t \geq 0; 
		\end{equation}
		 \begin{equation}\label{jhgtn0111} 
		 \int_{r}^{t}\|u(s)\|^{2\alpha+2}_{L^{2\alpha+2}}ds \leq \frac{1}{a}\|u(r)\|^{2}_{L^{2}} + \left(\frac{2l^{2}}{a\nu}\|f\|^{2}_{L^{\infty}_{t}L^{2}} +
		 \frac{2^{\frac{1+\alpha}{\alpha}}\nu^{\frac{\alpha+1}{\alpha}} l^{\frac{\alpha-2}{\alpha}}}{a^{\frac{1+\alpha}{\alpha}}}\right)(t-r), \,\, \forall\, t  \geq r \geq 0;
		 \end{equation}
	 
	 \begin{eqnarray}\label{jhgtn03} 
	 	\|\nabla u(t)\|_{L^{2}}^{2}  &\leq \left(\frac{1}{2l^{2}}+\frac{1}{2\nu^{\frac{2\alpha-1}{\alpha-1}} a^{\frac{1}{\alpha-1}}}\right)e^{\frac{-2\nu \left(t-\frac{l^{2}}{\nu}\right)}{l^{2}}}\|u(0)\|^{2}_{L^{2}} + K \left( \frac{3}{2\nu} +\frac{3l^{2}}{2\nu^{\frac{3\alpha-2}{\alpha-1}} a^{\frac{1}{\alpha-1}}} \right) \nonumber\\&+ \ds\frac{l^{2}}{\nu^{2}}\|f\|^{2}_{L^{\infty}_{t}L^{2}} ,  \,\, \forall\, t \geq \frac{l^{2}}{\nu};
	 \end{eqnarray}
		
	\begin{eqnarray}\label{jhgtn03.2}
	\|\nabla u(t)\|_{L^{2}}^{2}&\leq \|\nabla u(0)\|_{L^{2}}^{2} +\frac{1}{\nu^{\frac{\alpha}{\alpha-1}} a^{\frac{1}{\alpha-1}}}\left(\frac{l^{4}}{\nu^{3}}\|f\|^{2}_{L^{\infty}_{t}L^{2}}+ \frac{4 l^{\frac{3\alpha-2}{\alpha}}}{\nu^{\frac{\alpha-1}{\alpha}}  a^{\frac{1}{\alpha}}}\right) \nonumber\\&+\ds\frac{l^{2}}{\nu^{2}}\|f\|^{2}_{L^{\infty}_{t}L^{2}},\,\,   \forall\, 0\leq t \leq \frac{l^{2}}{\nu};
	\end{eqnarray}
		 
		\begin{equation} \label{jhgtn04}
		\int_{r}^{t}\|A u(s)\|_{L^{2}}^{2}ds\leq \frac{2}{\nu}\|\nabla u(r)\|_{L^{2}}^{2} +\frac{2}{\nu^{\frac{2\alpha-1}{\alpha-1}} a^{\frac{1}{\alpha-1}}} \int_{r}^{t}\|\nabla u(s)\|^{2}_{L^{2}}ds +\frac{4(t-r)}{\nu^{2}}\|f\|^{2}_{L^{\infty}_{t}L^{2}}, \,\, \forall\, t  \geq r \geq 0.	
		\end{equation}

\end{lee}


\subsubsection{Proof of Lemma \ref{da0101}}


Multiplying the system (\ref{eq111}) by the strong solution $u(t)$, integrating in $\Omega$, performing integration by parts and using (\ref{fdr2}), we  obtain the following equality:
\begin{equation}
	\frac{1}{2}\frac{d}{dt} \|u\|_{L^{2}}^{2}+\nu\|\nabla u\|_{L^{2}}^{2}+ a\|u\|^{2\alpha+2}_{L^{2\alpha+2}}= (f,u)_{L^{2}}.   \label{in1}
\end{equation}%
Using Young and Hölder's inequalities, we obtain
\begin{align}
	|(f,u)_{L^{2}}| &\leq \|u\|_{L^{2\alpha+2}}\|f\|_{L^{\frac{2\alpha+2}{2\alpha+1}}} \leq \frac{a}{2}\|u\|^{2\alpha+2}_{L^{2\alpha+2}}+\left(\frac{2}{a}\right)^{\frac{1}{2\alpha+1}}\|f\|^{\frac{2\alpha+2}{2\alpha+1}}_{L^{\frac{2\alpha+2}{2\alpha+1}}} \nonumber\\&\leq \frac{a}{2}\|u\|^{2\alpha+2}_{L^{2\alpha+2}}+\left(\frac{2}{a}\right)^{\frac{1}{2\alpha+1}}\|f\|^{\frac{2\alpha+2}{2\alpha+1}}_{L^{2}}l^{\frac{3\alpha}{2\alpha+1}} \leq
	\frac{a}{2}\|u\|^{2\alpha+2}_{L^{2\alpha+2}}+\frac{l^{2}}{\nu}\|f\|^{2}_{L^{2}}+ \frac{2^{\frac{1}{\alpha}}\nu^{\frac{\alpha+1}{\alpha}}}{a^{\frac{1}{\alpha}} l^{\frac{2\alpha+2}{\alpha}}}l^{3}. \label{rthuf}
\end{align}
Then, plugging \eqref{rthuf} into \eqref{in1}, we have
\begin{align}
	\frac{1}{2}\frac{d}{dt} \|u\|_{L^{2}}^{2}+\nu\|\nabla u\|_{L^{2}}^{2}+ \frac{a}{2}\|u\|^{2\alpha+2}_{L^{2\alpha+2}} &\leq  \frac{l^{2}}{\nu}\|f\|^{2}_{L^{2}}+ \frac{2^{\frac{1}{\alpha}}\nu^{\frac{\alpha+1}{\alpha}}}{a^{\frac{1}{\alpha}} l^{\frac{2\alpha+2}{\alpha}}}l^{3}. \label{454fghtbflk}
\end{align}
Furthermore, note that
\begin{align*}
	\frac{\nu}{l^{2}} \|u\|^{2}_{L^{2}} &\leq \frac{\nu}{l^{2}} \|u\|^{2}_{L^{2\alpha+2}}\cdot l^{\frac{3\alpha}{\alpha+1}}  \leq \frac{a}{2}\|u\|^{2\alpha+2}_{L^{2\alpha+2}}+  \frac{2^{\frac{1}{\alpha}}\nu^{\frac{\alpha+1}{\alpha}}}{a^{\frac{1}{\alpha}}l^{\frac{2\alpha+2}{\alpha}} }l^{3}. 
\end{align*}
Therefore 
\begin{align}
	\frac{1}{2}\frac{d}{dt} \|u\|_{L^{2}}^{2}+\nu\|\nabla u\|_{L^{2}}^{2}+ \frac{\nu}{l^{2}} \|u\|^{2}_{L^{2}} &\leq  \frac{l^{2}}{\nu}\|f\|^{2}_{L^{2}}+ 2\frac{2^{\frac{1}{\alpha}}\nu^{\frac{\alpha+1}{\alpha}}}{a^{\frac{1}{\alpha}}l^{\frac{2\alpha+2}{\alpha}} }l^{3}. \label{zzzwe01}
\end{align}
Using Gronwall's inequality, we obtain
\begin{align}
	\|u(t)\|_{L^{2}}^{2} &\leq e^{-\frac{2\nu t}{l^{2}}}\|u(0)\|_{L^{2}}^{2} + \frac{l^{2}}{\nu}\left(\frac{l^{2}}{\nu}\|f\|^{2}_{L^{\infty}_{t}L^{2}}+ 2\frac{2^{\frac{1}{\alpha}}\nu^{\frac{\alpha+1}{\alpha}}}{a^{\frac{1}{\alpha}}l^{\frac{2\alpha+2}{\alpha}} }l^{3}\right)\nonumber\\&\leq  e^{-\frac{2\nu t}{l^{2}}}\|u(0)\|^{2}_{L^{2}}+  \frac{l^{2}}{\nu}K, 
	\nonumber
\end{align}
concluding (\ref{jhgtn01}). Besides, if we integrate (\ref{454fghtbflk}) over $[r,t]$, we obtain (\ref{jhgtn0111}).

Now, integrating (\ref{zzzwe01}) over $[c,d]$ and using (\ref{jhgtn01}), we obtain
\begin{align}
	\int_{c}^{d}\|\nabla u(s)\|_{L^{2}}^{2}ds&\leq  \frac{1}{2\nu}\|u(c)\|^{2}_{L^{2}}+ (d-c)\frac{l^{2}}{\nu^{2}}\|f\|^{2}_{L^{\infty}_{t}L^{2}}+ (d-c)2\frac{2^{\frac{1}{\alpha}}\nu^{\frac{\alpha+1}{\alpha}-1}}{a^{\frac{1}{\alpha}}l^{\frac{2\alpha+2}{\alpha}} }l^{3}\nonumber\\&\leq \frac{1}{\nu}\left(\frac{1}{2}\|u(c)\|^{2}_{L^{2}}+ (d-c)K\right) \nonumber\\&\leq \frac{1}{\nu}\left(\frac{1}{2}\left(e^{-\frac{2\nu c}{l^{2}}}\|u(0)\|^{2}_{L^{2}}+  \frac{l^{2}}{\nu}K\right)+ (d-c)K\right). \label{poijn001}
\end{align}
Considering $c=t$ and $d =t+\frac{l^{2}}{\nu}$, we have
\begin{align}
	\int_{t}^{t+\frac{l^{2}}{\nu}}\|\nabla u(s)\|_{L^{2}}^{2}ds&\leq  \frac{1}{2\nu}e^{-\frac{2\nu t}{l^{2}}}\|u(0)\|^{2}_{L^{2}}+ \frac{3l^{2}}{2\nu^{2}}K,\,\, \forall\, t \geq 0. \label{jhgtn02}
\end{align}

Multiplying the system (\ref{eq111}) by $A u$, integrating over $\Omega$ and  performing integration by parts,  we  obtain the following equality: 
\begin{align}
	\frac{1}{2}\frac{d}{dt} \|\nabla u\|_{L^{2}}^{2}&+\nu\|A u\|_{L^{2}}^{2}+ a\||u|^{\alpha}|\nabla u|\|^{2}_{L^{2}}+\frac{a \alpha}{2}\||u|^{\alpha-1}\nabla |u|^{2}\|^{2}_{L^{2}} = -(u\cdot \nabla u, A u)_{L^{2}}+ (f,Au)_{L^{2}}.  
	\nonumber
\end{align}%
Using Young and Hölder's inequalities of the right side above, we get for $\epsilon>0$,
\begin{align}
	|(f,Au)_{L^{2}}| &\leq \|f\|_{L^{2}}\|Au\|_{L^{2}} \leq \nu\epsilon\|Au\|_{L^{2}}^{2}+\frac{1}{4\epsilon \nu}\|f\|^{2}_{L^{2}}, \label{rthuf3}\\\nonumber&\\
	|(u\cdot \nabla u, A u)_{L^{2}}|&\leq \int_{\Omega} |u| |\nabla u|^{\frac{1}{\alpha}} |\nabla u|^{\frac{\alpha-1}{\alpha}} |A u|dx \leq \||u||\nabla u|^{\frac{1}{\alpha}}\|_{L^{2\alpha}}  \||\nabla u|^{\frac{\alpha-1}{\alpha}} \|_{L^{\frac{2\alpha}{\alpha-1}}}  \|A u\|_{L^{2}} \nonumber\\&\leq \frac{1}{2\nu}\||u|^{\alpha}|\nabla u|\|^{\frac{2}{\alpha}}_{L^{2}}  \|\nabla u\|^{\frac{2(\alpha-1)}{\alpha}}_{L^{2}} + \frac{\nu}{2}\|A u\|^{2}_{L^{2}}  \nonumber\\&\leq \frac{a}{2}\||u|^{\alpha}|\nabla u|\|^{2}_{L^{2}}  +\left(\frac{1}{2\nu}\right)^{\frac{\alpha}{\alpha-1}}\left(\frac{2}{a}\right)^{\frac{1}{\alpha-1}}\|\nabla u\|^{2}_{L^{2}} + \frac{\nu}{2}\|A u\|^{2}_{L^{2}}; \label{rthuf4}
\end{align}
Hence we have
\begin{align}
	\frac{1}{2}\frac{d}{dt} \|\nabla u\|_{L^{2}}^{2}&+\nu\left(\frac{1}{2}-\epsilon\right)\|A u\|_{L^{2}}^{2}+ \frac{a}{2}\||u|^{\alpha}|\nabla u|\|^{2}_{L^{2}}+\frac{a\alpha}{2}\||u|^{\alpha-1}\nabla |u|^{2}\|^{2}_{L^{2}} \nonumber\\&\leq  \frac{1}{2\nu^{\frac{\alpha}{\alpha-1}} a^{\frac{1}{\alpha-1}}} \|\nabla u\|^{2}_{L^{2}}+\frac{1}{4\epsilon \nu}\|f\|^{2}_{L^{2}}.  \label{rgn1.1.1}
\end{align}
Fix $t \geq \frac{l^{2}}{\nu}$ and $r$ such that  $t-\frac{l^{2}}{\nu}\leq r \leq t$. Choosing $\epsilon=\frac{1}{2}$, integrating (\ref{rgn1.1.1}) over $[r,t]$ and using (\ref{poijn001}) properly,  we have
\begin{align}
	\frac{1}{2}\|\nabla u(t)\|_{L^{2}}^{2} &\leq \frac{1}{2}\|\nabla u(r)\|_{L^{2}}^{2} + \frac{1}{2\nu^{\frac{\alpha}{\alpha-1}} a^{\frac{1}{\alpha-1}}}\int_{t-\frac{l^{2}}{\nu}}^{t}\|\nabla u(s)\|^{2}_{L^{2}}ds +\frac{1}{2 \nu}\int_{t-\frac{l^{2}}{\nu}}^{t}\|f(s)\|^{2}_{L^{2}}ds \nonumber \\&\leq \frac{1}{2}\|\nabla u(r)\|_{L^{2}}^{2} + \frac{l^{2}}{2\nu^{2}}\|f\|^{2}_{L^{\infty}_{t}L^{2}}+ \frac{1}{2\nu^{\frac{\alpha}{\alpha-1}} a^{\frac{1}{\alpha-1}}} \left(\frac{1}{2\nu}e^{\frac{-2\nu (t-\frac{l^{2}}{\nu})}{l^{2}}}\|u(0)\|^{2}_{L^{2}}+  \frac{3l^{2}}{2\nu^{2}}K \right). \nonumber
\end{align}%
Now, integrating the above inequality with respect to $r$ over $[t-\frac{l^{2}}{\nu},t]$ and using again (\ref{poijn001}) appropriately, we obtain 
\begin{align}
	\frac{l^{2}}{\nu}\frac{1}{2}\|\nabla u(t)\|_{L^{2}}^{2}  &\leq \frac{1}{2}\int^{t}_{t-\frac{l^{2}}{\nu}}\|\nabla u(r)\|_{L^{2}}^{2}dr +\frac{l^{4}}{2\nu^{3}}\|f\|^{2}_{L^{\infty}_{t}L^{2}} +
	\frac{l^{2}}{2\nu^{\frac{3\alpha-2}{\alpha-1}} a^{\frac{1}{\alpha-1}}}\left(\frac{1}{2}e^{\frac{-2\nu \left(t-\frac{l^{2}}{\nu}\right)}{l^{2}}}\|u(0)\|^{2}_{L^{2}}+  \frac{3l^{2}}{2\nu}K \right) 
	\nonumber \\&\leq\frac{1}{2\nu}\left(\frac{1}{2}e^{\frac{-2\nu \left(t-\frac{l^{2}}{\nu}\right)}{l^{2}}}\|u(0)\|^{2}_{L^{2}}+ \frac{3l^{2}}{2\nu}K\right) + \frac{l^{4}}{2\nu^{3}}\|f\|^{2}_{L^{\infty}_{t}L^{2}}  \nonumber\\ &+
	\frac{l^{2}}{2\nu^{\frac{3\alpha-2}{\alpha-1}} a^{\frac{1}{\alpha-1}}}\left(\frac{1}{2}e^{\frac{-2\nu \left(t-\frac{l^{2}}{\nu}\right)}{l^{2}}}\|u(0)\|^{2}_{L^{2}}+  \frac{3l^{2}}{2\nu}K \right). \nonumber
\end{align}%
Therefore, we get (\ref{jhgtn03}). Moreover, integrating (\ref{zzzwe01}) over $[0, \frac{l^{2}}{\nu}]$, we obtain
\begin{align}
	\int_{0}^{\frac{l^{2}}{\nu}}\|\nabla u(s)\|_{L^{2}}^{2}ds&\leq  \frac{l^{4}}{\nu^{3}}\|f\|^{2}_{L^{\infty}_{t}L^{2}}+ \frac{4 l^{\frac{3\alpha-2}{\alpha}}}{\nu^{\frac{\alpha-1}{\alpha}}  a^{\frac{1}{\alpha}}}. \label{zzzwe01.dfsrvrvsd}
\end{align}
Let  $0 \leq t \leq \frac{l^{2}}{\nu}$. Integrating (\ref{rgn1.1.1}) with $\epsilon=\frac{1}{2}$ over $[0, t]$, we obtain
\begin{align}
	\|\nabla u(t)\|_{L^{2}}^{2}&\leq \|\nabla u(0)\|_{L^{2}}^{2} +\frac{1}{\nu^{\frac{\alpha}{\alpha-1}} a^{\frac{1}{\alpha-1}}}\left(\frac{l^{4}}{\nu^{3}}\|f\|^{2}_{L^{\infty}_{t}L^{2}}+ \frac{4 l^{\frac{3\alpha-2}{\alpha}}}{\nu^{\frac{\alpha-1}{\alpha}}  a^{\frac{1}{\alpha}}}\right) +\frac{l^{2}}{\nu^{2}}\|f\|^{2}_{L^{\infty}_{t}L^{2}}.  \label{rgn1.1.1kjkjlkj}
\end{align}
and thus (\ref{jhgtn03.2}) is valid. 

Finally, integrating (\ref{rgn1.1.1}) over $[r,t]$ choosing $\epsilon = \frac{1}{4}$ we obtain (\ref{jhgtn04}).

\begin{flushright}
	\rule{2mm}{2mm}
\end{flushright}

We summarize Lemma \ref{da0101} below:

\begin{col} \label{da0101.4565.1}
	Suppose $f\in L^{\infty}(\mathbb{R}_{+};H)$, $\alpha>1$ and  $a>0$. Let
	$u$ be a global strong solution of (\ref{eq111}). Consider $K$  as given in (\ref{ct01}) and $M >0$ satisfying (\ref{estM}). Define the following
		$$M_{1}= l^{2}M + \frac{l^{2}K}{\nu},$$
	$$	\tilde{M}=\max\left\{\frac{1}{\nu^{\frac{\alpha}{\alpha-1}} a^{\frac{1}{\alpha-1}}}\left(\frac{l^{4}}{\nu^{3}}\|f\|^{2}_{L^{\infty}_{t}L^{2}}+ \frac{4 l^{\frac{3\alpha-2}{\alpha}}}{\nu^{\frac{\alpha-1}{\alpha}}a^{\frac{1}{\alpha}}}  \right), \,\, \frac{l^{2}}{2\nu^{\frac{2\alpha-1}{\alpha-1}} a^{\frac{1}{\alpha-1}}}M + K \left( \frac{3}{2\nu} +\frac{3l^{2}}{2\nu^{\frac{3\alpha-2}{\alpha-1}} a^{\frac{1}{\alpha-1}}} \right) \right\}, $$
	$$	M_{2}= \frac{l^{2}}{\nu^{2}}\|f\|^{2}_{L^{\infty}_{t}L^{2}} + M+\tilde{M},$$
	$$	M_{3}= \frac{2}{\nu^{\frac{2\alpha-1}{\alpha-1}} a^{\frac{1}{\alpha-1}}}M_{2} +\frac{4}{\nu^{2}}\|f\|^{2}_{L^{\infty}_{t}L^{2}},$$
	$$	M_{4}= \frac{2l^{2}}{a\nu}\|f\|^{2}_{L^{\infty}_{t}L^{2}} +
		\frac{2^{\frac{1+\alpha}{\alpha}}\nu^{\frac{\alpha+1}{\alpha}} l^{\frac{\alpha-2}{\alpha}}}{a^{\frac{1+\alpha}{\alpha}}}.$$
	
	Then,
	\begin{align}
		\|u(t)\|_{L^{2}}^{2} &\leq M_{1}, \,\, \forall\, t \geq 0, \nonumber\\&\nonumber\\
		\|\nabla u(t)\|_{L^{2}}^{2}  &\leq  M_{2} ,  \,\, \forall\, t \geq 0, \nonumber\\&\nonumber\\
		\int_{r}^{t}\|A u(s)\|_{L^{2}}^{2}ds&\leq \frac{2}{\nu}M_{2}+ (t-r)M_{3}, \,\, \forall\, t \geq r  \geq 0, \nonumber\\&\nonumber\\
		\int_{r}^{t}\|u(s)\|^{2\alpha+2}_{L^{2\alpha+2}}ds &\leq \frac{1}{a}M_{1} + (t-r)M_{4}, \,\, \forall\, t  \geq r \geq 0. \nonumber
	\end{align}
\end{col}


\begin{col} \label{da0101.4565.1.2}
	Suppose $f\in L^{\infty}(\mathbb{R}_{+};H)$, $\alpha>1$,  $a>0$ and 
	$u$ be a global strong solution of (\ref{eq111}). Consider $M >0$ that satisfies (\ref{estM})
	and $M_{2}$ as in Corollary \ref{da0101.4565.1}. 
	Then, fixed $\eta>0$, we have for all $t\geq 0$,
	\begin{align}
		\int_{0}^{t}e^{\frac{\eta}{8}(s-t)}\|A u(s)\|_{L^{2}}^{2}ds\leq \frac{4}{\nu}M_{2}+ \frac{16}{\eta \nu^{\frac{3\alpha-2}{\alpha-1}} a^{\frac{1}{\alpha-1}}} M_{2} +\frac{32}{\eta \nu^{2}}\|f\|^{2}_{L^{\infty}_{t}L^{2}},  \nonumber
	\end{align}
\end{col}


\subsubsection{Proof of Corollary \ref{da0101.4565.1.2}}


Using inequality (\ref{rgn1.1.1}) with $\epsilon= \frac{1}{4}$ we obtain
\begin{align}
	\frac{d}{ds} \|\nabla u(s)\|_{L^{2}}^{2}&+\frac{\nu}{2}\|A u(s)\|_{L^{2}}^{2}\leq  \frac{1}{\nu^{\frac{\alpha}{\alpha-1}} a^{\frac{1}{\alpha-1}}} M_{2}+\frac{2}{\nu}\|f\|^{2}_{L^{\infty}_{t}L^{2}}. \nonumber 
\end{align} 
Multiply the above inequality by $e^{\frac{\eta}{8}(s-t)}$, we have
\begin{align}
	\frac{d}{ds} \left(e^{\frac{\eta}{8}(s-t)}\|\nabla u(s)\|_{L^{2}}^{2}\right)+\frac{\nu}{2}e^{\frac{\eta}{8}(s-t)}\|A u(s)\|_{L^{2}}^{2}\leq  e^{\frac{\eta}{8}(s-t)}\left(\frac{1}{\nu^{\frac{\alpha}{\alpha-1}} a^{\frac{1}{\alpha-1}}} M_{2} +\frac{\eta}{8}M_{2}+\frac{2}{\nu}\|f\|^{2}_{L^{\infty}_{t}L^{2}}\right). \nonumber 
\end{align} 
Integrating over $[0,t]$, we obtain
\begin{align}
	\|\nabla u(t)\|_{L^{2}}^{2} +\frac{\nu}{2}\int_{0}^{t}e^{\frac{\eta}{8}(s-t)}\|A u(s)\|_{L^{2}}^{2}ds\leq e^{\frac{-\eta}{8}t}\|\nabla u(0)\|_{L^{2}}^{2}+ \frac{8}{\eta}\left(\frac{1}{\nu^{\frac{\alpha}{\alpha-1}} a^{\frac{1}{\alpha-1}}} M_{2} +\frac{\eta}{8}M_{2}+\frac{2}{\nu}\|f\|^{2}_{L^{\infty}_{t}L^{2}}\right).  \nonumber
\end{align} 
Then, we have the result.

\begin{flushright}
	\rule{2mm}{2mm}
\end{flushright}


\begin{lee} \label{da0101.AA}
	Suppose $f\in L^{\infty}(\mathbb{R}_{+};H)$ with $f_{t}\in L^{\infty}(\mathbb{R}_{+};H)$, $1 <  \alpha < 2$, $a>0$ and 
	$u$ be a global strong solution of (\ref{eq111}). Consider $M >0$ such that (\ref{estM}) is satisfied and $M_{1}$, $M_{2}$, $M_{3}$ and $M_{4}$ as given in Corollary \ref{da0101.4565.1}. Let 
	\begin{align}
		M_{5}&= 
		\frac{\nu}{al^{2}}M_{1}+ M_{4}+  \frac{\nu(\alpha+1)}{a}M_{2}+ \frac{(2\alpha+2)l^{2}}{a \nu}\|f\|^{2}_{L^{\infty}_{t}L^{2}} \nonumber\\&+
		\frac{4C_{\infty}^{2} (2\alpha+2)}{a}\left(\frac{4 l^{2}}{\nu^{3}}M^{3}_{2}+2\nu M_{2}+l^{2}\nu M_{3}+ \frac{1}{l^{2}\nu}M_{1}M_{2}\right), \nonumber \\
		M_{6}&=\nu M_{2}+  \frac{a}{\alpha+1}M_{5}+ \frac{2 l^{2}}{\nu} \|f\|^{2}_{L^{\infty}_{t}L^{2}} + 8C_{\infty}^{2}\left(\frac{4l^{2}}{\nu^{3}}M^{3}_{2}+2\nu M_{2}+l^{2}\nu M_{3}+ \frac{1}{l \nu}M_{1}M_{2}\right), \nonumber\\
		M_{7}&= M_{6}\left[\frac{3\nu}{2l^{2}}+ \frac{108(C_{3}C_{6})^{4}}{\nu^{3}}M_{2}^{2}+\frac{6(C_{3}C_{6})^{\frac{4}{3}}}{\nu^{\frac{1}{3}} l^{\frac{4}{3}}}M_{2}^{\frac{2}{3}} +\frac{4C^{2}_{3}C^{2}_{6}}{\nu l}M_{2}+ \frac{2C_{3}C_{6}}{ l^{\frac{3}{2}}}M_{2}^{\frac{1}{2}}\right]+\frac{2l^{4}}{\nu^{2}}\|f_{t}\|_{L^{\infty}_{t}L^{2}}^{2}, \nonumber\\
	\end{align}
\begin{eqnarray}
		M_{8}= \frac{2}{\nu}M_{7}+ 2^{10}\frac{C^{8}_{4}}{\nu^{4}} \big(M^{\frac{1}{2}}_{1} M_{2}^{\frac{3}{2}} &\!\!\!\!\!+&\!\!\!\!\!\frac{1}{l^{3}}M_{1}^{2}\big)M^{\frac{1}{2}}_{2}+ \frac{2^{\frac{4\alpha}{2-\alpha}}a^{\frac{2}{2-\alpha}}M^{\frac{1}{2-\alpha}}_{5}C^{\frac{2\alpha}{2-\alpha}}_{\infty}}{\nu^{\frac{2}{2-\alpha}}}M_{2}^{\frac{\alpha}{4-2\alpha}}\nonumber\\
		 &\!\!\!\!\!+&\!\!\!\!\! \frac{2^{\alpha}aM^{\frac{1}{2}}_{5}C^{\alpha}_{\infty}}{\nu l^{\frac{3\alpha}{2}}}M_{1}^{\frac{\alpha}{2}} +\frac{2}{\nu}\|f\|_{L^{\infty}_{t}L^{2}}. \nonumber
	\end{eqnarray}
	Then,
	\begin{align}
		\left\|u(t)\right\|^{2\alpha+2}_{L^{2\alpha+2}} &\leq M_{5},\,\, \forall\, t \geq \frac{l^{2}}{\nu}; \label{jhgtn02.AAZZ} \\& \nonumber\\
		\int_{t}^{t+\frac{l^{2}}{\nu}}\|u_{t}(s)\|^{2}_{L^{2}} ds &\leq M_{6},\,\, \forall\, t \geq \frac{l^{2}}{\nu}; \label{jhgtn02.AAZZ11} \\& \nonumber\\
		\| u_{t}(t)\|_{L^{2}}^{2} &\leq M_{7}, \,\, \forall\, t \geq \frac{2l^{2}}{\nu}; \label{jhgtn02.AAZZ.2}\\& \nonumber\\
		\| A u(t)\|_{L^{2}} &\leq M_{8},  \,\, \forall\, t \geq \frac{2l^{2}}{\nu}. \label{jhgtn02.AAZZ.3}
	\end{align}
\end{lee}


\subsubsection{Proof of Lemma \ref{da0101.AA}}


Multiplying the system (\ref{eq111}) by $u_{t}$, integrating in $\Omega$ and  performing integration by parts, we obtain the following: 
\begin{align}
	\| u_{t}\|_{L^{2}}^{2}+\frac{\nu}{2}\frac{d}{dt}\|\nabla u\|_{L^{2}}^{2}+ \frac{a}{(2\alpha+2)}\frac{d}{dt}\||u|^{\alpha+1}\|^{2}_{L^{2}}= -(u\cdot \nabla u,u_{t})_{L^{2}}+(f,u_{t})_{L^{2}}.   \label{in2.5}
\end{align}%
By (\ref{gn}) and Hölder's inequality, we have
\begin{align}
	|(u\cdot \nabla u, u_{t})_{L^{2}}|&\leq \|u_{t}\|_{L^{2}}\|u\|_{L^{\infty}}\|\nabla u\|_{L^{2}} \leq C_{\infty}\|u_{t}\|_{L^{2}}\left(\|\nabla u\|^{\frac{1}{2}}_{L^{2}}\|A u\|^{\frac{1}{2}}_{L^{2}}+ \frac{1}{l^{\frac{3}{2}}}\|u\|_{L^{2}}\right)\|\nabla u\|_{L^{2}}\nonumber\\ &\leq \epsilon \|u_{t}\|^{2}_{L^{2}}+ \frac{C_{\infty}^{2}}{\epsilon}\left(\|\nabla u\|^{3}_{L^{2}}\|A u\|_{L^{2}}+ \frac{1}{l^{3}}\|u\|^{2}_{L^{2}}\|\nabla u\|^{2}_{L^{2}}\right), \nonumber\\&\nonumber
\end{align}
and
\begin{equation*}
		|(f,u_{t})_{L^{2}}| \leq  \epsilon\|u_{t}\|_{L^{2}}^{2}+\frac{1}{4\epsilon}\|f\|^{2}_{L^{2}}.
\end{equation*}
Thus,
\begin{align}
	(1-2\epsilon)\| u_{t}\|_{L^{2}}^{2}+\frac{\nu}{2}\frac{d}{dt}\|\nabla u\|_{L^{2}}^{2}&+ \frac{a}{(2\alpha+2)}\frac{d}{dt}\||u|^{\alpha+1}\|^{2}_{L^{2}} \leq    \frac{1}{4\epsilon}\|f\|^{2}_{L^{2}} \nonumber\\&+ \frac{C_{\infty}^{2}}{\epsilon}\left(\|\nabla u\|^{3}_{L^{2}}\|A u\|_{L^{2}}+ \frac{1}{l^{3}}\|u\|^{2}_{L^{2}}\|\nabla u\|^{2}_{L^{2}}\right). \nonumber 
\end{align}
Choosing $\epsilon=\frac{1}{4}$, Hölder's inequality implies
\begin{align}
	\frac{1}{2} \| u_{t}\|_{L^{2}}^{2}+\frac{\nu}{2}\frac{d}{dt}\|\nabla u\|_{L^{2}}^{2}&+ \frac{a}{(2\alpha+2)}\frac{d}{dt}\|u\|^{2\alpha+2}_{L^{2\alpha+2}} \leq    \|f\|^{2}_{L^{2}} \nonumber\\&+ 4C_{\infty}^{2}\left(\frac{4}{\nu^{2}}\|\nabla u\|^{6}_{L^{2}}+\nu^{2}\|A u\|^{2}_{L^{2}}+ \frac{1}{l^{3}}\|u\|^{2}_{L^{2}}\|\nabla u\|^{2}_{L^{2}}\right).  \label{iiii6} 
\end{align}
Fix now $t \geq  0$. Considering $s$ such that $t< s<t+\frac{l^{2}}{\nu}$, integrating over $\big[s,t+\frac{l^{2}}{\nu}\big]$ and using Corollary  \ref{da0101.4565.1}, we obtain
\begin{align}
	\frac{a}{(2\alpha+2)}\left\|u\left(t+\frac{l^{2}}{\nu}\right)\right\|^{2\alpha+2}_{L^{2\alpha+2}} &\leq
	\frac{a}{(2\alpha+2)}\|u(s)\|^{2\alpha+2}_{L^{2\alpha+2}}+ \frac{\nu}{2}\|\nabla u(s)\|_{L^{2}}^{2} + \frac{l^{2}}{\nu}\|f\|^{2}_{L^{\infty}_{t}L^{2}} \nonumber\\&+
	4C_{\infty}^{2}\int_{t}^{t+\frac{l^{2}}{\nu}}\left(\frac{4}{\nu^{2}}\|\nabla u(r)\|^{6}_{L^{2}}+\nu^{2}\|A u(r)\|^{2}_{L^{2}}+ \frac{1}{l^{3}}\|u(r)\|^{2}_{L^{2}}\|\nabla u(r)\|^{2}_{L^{2}}\right)dr\nonumber\\&\leq 
	\frac{a}{(2\alpha+2)}\|u(s)\|^{2\alpha+2}_{L^{2\alpha+2}}+ \frac{\nu}{2}M_{2} + \frac{l^{2}}{\nu}\|f\|^{2}_{L^{\infty}_{t}L^{2}} \nonumber\\&+
	4C_{\infty}^{2}\left(\frac{4 l^{2}}{\nu^{3}}M^{3}_{2}+2\nu M_{2}+l^{2}\nu M_{3}+ \frac{1}{l\nu}M_{1}M_{2}\right). \nonumber
\end{align}%
Integrating above over $\big[t,t+\frac{l^{2}}{\nu}\big]$ in $s$, we get
\begin{align}
	\frac{a l^{2}}{\nu(2\alpha+2)}\left\|u\left(t+\frac{l^{2}}{\nu}\right)\right\|^{2\alpha+2}_{L^{2\alpha+2}} &\leq 
	\frac{a}{(2\alpha+2)}\int_{t}^{t+\frac{l^{2}}{\nu}}\|u(s)\|^{2\alpha+2}_{L^{2\alpha+2}}ds+  \frac{l^{2}}{2}M_{2} + \frac{l^{4}}{\nu^{2}}\|f\|^{2}_{L^{\infty}_{t}L^{2}} \nonumber\\&+
	4C_{\infty}^{2}\left(\frac{4 l^{4}}{\nu^{4}}M^{3}_{2}+2l^{2} M_{2}+l^{4}M_{3}+ \frac{l}{\nu^{2}}M_{1}M_{2}\right)
	\nonumber\\&\leq 
	\frac{a}{(2\alpha+2)}\left(\frac{1}{a}M_{1}+ \frac{l^{2}}{\nu}M_{4}\right)+ \frac{l^{2}}{2}M_{2} + \frac{l^{4}}{\nu^{2}}\|f\|^{2}_{L^{\infty}_{t}L^{2}} \nonumber\\&+
	4C_{\infty}^{2}\left(\frac{4 l^{4}}{\nu^{4}}M^{3}_{2}+2l^{2} M_{2}+l^{4}M_{3}+ \frac{l}{\nu^{2}}M_{1}M_{2}\right).\nonumber
\end{align}%
Then, we obtain (\ref{jhgtn02.AAZZ}).

Consider $t \geq \frac{l^{2}}{\nu}$. Integrating (\ref{iiii6}) over $\big[t,t+\frac{l^{2}}{\nu}\big]$, we have
\begin{align}
	\frac{1}{2}\int_{t}^{t+\frac{l^{2}}{\nu}} \| u_{t}(s)\|_{L^{2}}^{2}ds&\leq \frac{\nu}{2}\|\nabla u(t)\|_{L^{2}}^{2}+  \frac{a}{2\alpha+2}\|u(t)\|^{2\alpha+2}_{L^{2\alpha+2}}+ \frac{l^{2}}{\nu} \|f\|^{2}_{L^{\infty}_{t}L^{2}} \nonumber\\&+ 4C_{\infty}^{2}\left(\frac{4l^{2}}{\nu^{3}}M^{3}_{2}+2\nu M_{2}+l^{2}\nu M_{3}+ \frac{1}{l \nu}M_{1}M_{2}\right)
	\nonumber \\ &\leq \frac{\nu}{2}M_{2}+  \frac{a}{2\alpha+2}M_{5}+ \frac{l^{2}}{\nu} \|f\|^{2}_{L^{\infty}_{t}L^{2}} \nonumber\\&+ 4C_{\infty}^{2}\left(\frac{4l^{2}}{\nu^{3}}M^{3}_{2}+2\nu M_{2}+l^{2}\nu M_{3}+ \frac{1}{l \nu}M_{1}M_{2}\right). \nonumber 
\end{align}
Therefore we conclude (\ref{jhgtn02.AAZZ11}).

We now differentiate (\ref{eq111}) with respect to $t$ and take the inner product with $u_{t}$ to obtain
\begin{align}
	\frac{1}{2}\frac{d}{dt}\| u_{t}\|_{L^{2}}^{2}+\frac{\nu}{2}\|\nabla u_{t}\|_{L^{2}}^{2}+ a(2\alpha+1)\||u|^{\alpha}u_{t}\|^{2}_{L^{2}} = -(u_{t}\cdot \nabla u,u_{t})_{L^{2}}+ (f_{t},u_{t})_{L^{2}}.   \label{457457ffd}
\end{align}%
Using (\ref{gn}) and Hölder's inequality, 
\begin{equation*}
	(f_{t},u_{t})_{L^{2}} \leq \frac{\nu}{4l^{2}}\| u_{t}\|_{L^{2}}^{2} + \frac{l^{2}}{\nu}\|f_{t}\|_{L^{2}}^{2},
	\end{equation*}
and
	\begin{align}
	\left|(u_{t}\cdot \nabla u,u_{t})_{L^{2}}\right|&\leq \|u_{t}\|_{L^{3}}\|\nabla u\|_{L^{2}}\|u_{t}\|_{L^{6}} \nonumber\\&\leq C_{3}C_{6}\left(\|u_{t}\|^{\frac{1}{2}}_{L^{2}}\|\nabla u_{t}\|^{\frac{1}{2}}_{L^{2}} + \frac{1}{l^{\frac{1}{2}}}\|u_{t}\|_{L^{2}}\right)\|\nabla u\|_{L^{2}}\left( \|\nabla u_{t}\|_{L^{2}} +\frac{1}{l}\|u_{t}\|_{L^{2}}\right)
	\nonumber \\&\leq \left(\frac{54(C_{3}C_{6})^{4}}{\nu^{3}}\|\nabla u\|^{4}_{L^{2}}+\frac{3(C_{3}C_{6})^{\frac{4}{3}}}{\nu^{\frac{1}{3}} l^{\frac{4}{3}}}\|\nabla u\|^{\frac{4}{3}}_{L^{2}} +\frac{2C^{2}_{3}C^{2}_{6}}{\nu l}\|\nabla u\|^{2}_{L^{2}}+ \frac{C_{3}C_{6}}{ l^{\frac{3}{2}}}\|\nabla u\|_{L^{2}}\right)\|u_{t}\|^{2}_{L^{2}}\nonumber\\&+\frac{3\nu}{8}\|\nabla u_{t}\|_{L^{2}}^{2}. \nonumber
\end{align}
Then,  by (\ref{457457ffd}), 
\begin{equation}
	\frac{1}{2}\frac{d}{dt}\| u_{t}\|_{L^{2}}^{2} \leq \left(\frac{\nu}{4l^{2}}+ \frac{54(C_{3}C_{6})^{4}}{\nu^{3}}M_{2}^{2}+\frac{3(C_{3}C_{6})^{\frac{4}{3}}}{\nu^{\frac{1}{3}} l^{\frac{4}{3}}}M_{2}^{\frac{2}{3}} +\frac{2C^{2}_{3}C^{2}_{6}}{\nu l}M_{2}+ \frac{C_{3}C_{6}}{ l^{\frac{3}{2}}}M_{2}^{\frac{1}{2}}\right)\|u_{t}\|^{2}_{L^{2}}+\frac{l^{2}}{\nu}\|f_{t}\|_{L^{2}}^{2}. \nonumber
\end{equation}%
Fixed $t \geq \frac{l^{2}}{\nu}$, consider $s$ such that $t<s<t+\frac{l^{2}}{\nu}$. Integrating the above inequality over $\big[s,t+\frac{l^{2}}{\nu}\big]$, using (\ref{jhgtn02.AAZZ11}) and defining $$K_{2}= \frac{\nu}{4l^{2}}+ \frac{54(C_{3}C_{6})^{4}}{\nu^{3}}M_{2}^{2}+\frac{3(C_{3}C_{6})^{\frac{4}{3}}}{\nu^{\frac{1}{3}} l^{\frac{4}{3}}}M_{2}^{\frac{2}{3}} +\frac{2C^{2}_{3}C^{2}_{6}}{\nu l}M_{2}+ \frac{C_{3}C_{6}}{ l^{\frac{3}{2}}}M_{2}^{\frac{1}{2}},$$  we get
\begin{align}
	\left\| u_{t}\left(t+\frac{l^{2}}{\nu}\right)\right\|_{L^{2}}^{2} &\leq \| u_{t}(s)\|_{L^{2}}^{2}+ 2K_{2}\int_{t}^{t+\frac{l^{2}}{\nu}}\|u_{t}(r)\|^{2}_{L^{2}}dr +\frac{2l^{4}}{\nu^{2}}\|f_{t}\|_{L^{\infty}_{t}L^{2}}^{2} \nonumber\\&\leq \| u_{t}(s)\|_{L^{2}}^{2}+ 2K_{2} M_{6}+\frac{2l^{4}}{\nu^{2}}\|f_{t}\|_{L^{\infty}_{t}L^{2}}^{2}. \nonumber
\end{align}%
Integrating the above inequality over $\big[t,t+\frac{l^{2}}{\nu}\big]$ in $s$ and using again (\ref{jhgtn02.AAZZ11}), we obtain
\begin{eqnarray}
	\frac{l^{2}}{\nu}\left\| u_{t}\left(t+\frac{l^{2}}{\nu}\right)\right\|_{L^{2}}^{2} &\!\!\!\!\leq \!\!\!\!& \int_{t}^{t+\frac{l^{2}}{\nu}}\| u_{t}(s)\|_{L^{2}}^{2}ds +\frac{2l^{2}K_{2} M_{6}}{\nu}+ \frac{2l^{6}}{\nu^{3}}\|f_{t}\|_{L^{\infty}_{t}L^{2}}^{2}\nonumber\\
	 &\!\!\!\!\leq \!\!\!\!& M_{6}\left(1+\frac{2l^{2}K_{2}}{\nu}\right)+\frac{2l^{6}}{\nu^{3}}\|f_{t}\|_{L^{\infty}_{t}L^{2}}^{2}. \nonumber
\end{eqnarray}%
Then,  we conclude (\ref{jhgtn02.AAZZ.2}). Finally, we take the $L^{2}$-norm in (\ref{eq111})  and get
\begin{equation}
	\nu \|A u\|_{L^{2}} \leq \|u_{t}\|_{L^{2}}+\|(u\cdot\nabla)  u\|_{L^{2}}  +a\||u|^{2\alpha}u\|_{L^{2}} +\|f\|_{L^{2}}. \label{fdfdbhg}
\end{equation}
Using  (\ref{gn}) and Hölder and Young's inequalities, we have 
\begin{align}
	\|(u\cdot\nabla)  u\|_{L^{2}} &\leq \|u\|_{L^{4}}\|\nabla u\|_{L^{4}} \leq C^{2}_{4}\left(\|u\|^{\frac{1}{4}}_{L^{2}} \|\nabla u\|^{\frac{3}{4}}_{L^{2}} +\frac{1}{l^{\frac{3}{4}}}\|u\|_{L^{2}}\right)\|\nabla u\|^{\frac{1}{4}}_{L^{2}}\|A u\|^{\frac{3}{4}}_{L^{2}} \nonumber\\ &\leq 432\frac{C^{8}_{4}}{\nu^{3}} \left(\|u\|_{L^{2}} \|\nabla u\|^{3}_{L^{2}} +\frac{1}{l^{3}}\|u\|^{4}_{L^{2}}\right)\|\nabla u\|_{L^{2}} +\frac{\nu}{8}\|A u\|_{L^{2}}. \nonumber
\end{align}
Since $1 <\alpha< 2$, using (\ref{gn}) and (\ref{jhgtn02.AAZZ}), we obtain
\begin{align}
	a\||u|^{2\alpha}u\|_{L^{2}} &\leq a\||u|^{\alpha+1}\|_{L^{2}}\|u^{\alpha}\|_{L^{\infty}}\leq
	a\|u\|^{\alpha+1}_{L^{2\alpha+2}}C^{\alpha}_{\infty}\left(\|\nabla u\|^{\frac{1}{2}}_{L^{2}}\|A u\|^{\frac{1}{2}}_{L^{2}}+\frac{1}{l^{\frac{3}{2}}}\|u\|_{L^{2}}\right)^{\alpha}
	\nonumber\\&\leq 2^{\alpha -1}aM^{\frac{1}{2}}_{5}C^{\alpha}_{\infty}\left(\|\nabla u\|^{\frac{\alpha}{2}}_{L^{2}}\|A u\|^{\frac{\alpha}{2}}_{L^{2}}+\frac{1}{l^{\frac{3\alpha}{2}}}\|u\|^{\alpha}_{L^{2}}\right)
	\nonumber\\&\leq \frac{2^{\frac{5\alpha -2}{2-\alpha}}a^{\frac{2}{2-\alpha}}M^{\frac{1}{2-\alpha}}_{5}C^{\frac{2\alpha}{2-\alpha}}_{\infty}}{\nu^{\frac{\alpha}{2-\alpha}}}\|\nabla u\|^{\frac{\alpha}{2-\alpha}}_{L^{2}} + \frac{2^{\alpha -1}aM^{\frac{1}{2}}_{5}C^{\alpha}_{\infty}}{l^{\frac{3\alpha}{2}}}\|u\|^{\alpha}_{L^{2}} +\frac{\nu}{8}\|A u\|_{L^{2}}, \nonumber
\end{align}
for all $ t > \frac{l^{2}}{\nu}$. Next, using above estimates in (\ref{fdfdbhg}), we conclude that
\begin{align}
	\frac{6\nu}{8}\|A u\|_{L^{2}} &\leq M_{7}+ 432\frac{C^{8}_{4}}{\nu^{3}} \left(M^{\frac{1}{2}}_{1} M_{2}^{\frac{3}{2}} +\frac{1}{l^{3}}M_{1}^{2}\right)M^{\frac{1}{2}}_{2}+ \frac{2^{\frac{5\alpha -2}{2-\alpha}}a^{\frac{2}{2-\alpha}}M^{\frac{1}{2-\alpha}}_{5}C^{\frac{2\alpha}{2-\alpha}}_{\infty}}{\nu^{\frac{\alpha}{2-\alpha}}}M_{2}^{\frac{\alpha}{4-2\alpha}} \nonumber\\&+ \frac{2^{\alpha -1}aM^{\frac{1}{2}}_{5}C^{\alpha}_{\infty}}{l^{\frac{3\alpha}{2}}}M_{1}^{\frac{\alpha}{2}} +\|f\|_{L^{\infty}_{t}L^{2}}. \nonumber
\end{align}
Then, we have (\ref{jhgtn02.AAZZ.3}).

\begin{flushright}
	\rule{2mm}{2mm}
\end{flushright}


\section{ Proof of Theorem \ref{CDA3a1.b}}\label{sec5} 


Let $u$ be the strong solution to (\ref{eq111}) and $w$ the weak solution of (\ref{eq1112}). Denoting $g =w-u$, we have
\begin{align}
	\frac{dg}{dt}+\nu Ag+B(w,w)-B(u,u)+bG_{\beta}(w) - aG_{\alpha}(u)=-\eta\mathcal{P}(I_{h}(g)). \label{eq111.0} 
\end{align}
Multiplying the system (\ref{eq111.0}) by $g$, integrating in space, using integration by parts and (\ref{fdr3}), we obtain
\begin{align}
	\frac{1}{2}\frac{d}{dt} \|g\|_{L^{2}}^{2}&+\nu\|\nabla g\|_{L^{2}}^{2}+ b\left\langle|w|^{2\beta}w- |u|^{2\beta}u,g\right\rangle_{Y_{\beta}',Y_{\beta}} =  \tilde{a}\frac{\nu}{l^{2}}\left\langle\left|\frac{l}{\nu}u\right|^{2\alpha}u- \left|\frac{l}{\nu}u\right|^{2\beta}u,g\right\rangle_{Y_{\beta}',Y_{\beta}}\nonumber\\&+\left(\tilde{a}-\tilde{b}\right)\frac{\nu}{l^{2}}\left\langle \left|\frac{l}{\nu}u\right|^{2\beta}u,g\right\rangle_{Y_{\beta}',Y_{\beta}}+\frac{1}{2}\left\langle g\cdot \nabla g,w+u\right\rangle_{Y_{\beta}',Y_{\beta}}  -\eta(I_{h}(g)-g,g)_{L^{2}}-\eta\|g\|^{2}_{L^{2}}. \label{in11}
\end{align}

We now estimate the terms of (\ref{in11}). Using (\ref{fdr1})-(\ref{dam02}) along Young and Hölder's inequalities, we have the following estimates:
\begin{align}
	\left|\eta (I_{h}(g)-g,g)_{L^{2}}\right| \leq \eta\widehat{\epsilon}\|I_{h}(g)-g\|^{2}_{L^{2}} +  \eta\frac{1}{4\widehat{\epsilon}}\|g\|^{2}_{L^{2}}\leq  \eta c_{o}h^{2}\widehat{\epsilon}\|\nabla g\|_{L^{2}}^{2} + \eta\frac{1}{4\widehat{\epsilon}}\|g\|^{2}_{L^{2}}; \label{iii3}
\end{align}
\begin{align}
	\frac{1}{2}\left|\left\langle g\cdot\nabla g,w+u\right\rangle_{Y',Y}\right| &\leq  \frac{1}{2}\int_{\Omega} |g|^{\frac{1}{\beta}} |g|^{\frac{\beta-1}{\beta}}|u+w||\nabla g| dx \nonumber\\
	&\leq  \frac{1}{2}\||g|^{\frac{1}{\beta}}|u+w|\|_{L^{2\beta}} \||g|^{\frac{\beta-1}{\beta}}\|_{L^{\frac{2\beta}{\beta-1}}}\|\nabla g\|_{L^{2}}  
	=  \frac{1}{2}\||g||u+w|^{\beta}\|^{\frac{1}{\beta}}_{L^{2}} \|g\|^{\frac{\beta-1}{\beta}}_{L^{2}}\|\nabla g\|_{L^{2}}  \nonumber\\
	&\leq  \frac{1}{16\nu\epsilon}\||g||u+w|^{\beta}\|^{\frac{2}{\beta}}_{L^{2}} \|g\|^{2\frac{\beta-1}{\beta}}_{L^{2}}+\nu\epsilon \|\nabla g\|^{2}_{L^{2}}\nonumber\\
	&\leq  \frac{1}{\nu^{\beta}(\eta \check{\epsilon})^{\beta-1}}\||g||u+w|^{\beta}\|^{2}_{L^{2}}+ \frac{(\beta-1)\eta \check{\epsilon}}{\beta (16\epsilon)^{\frac{\beta}{\beta-1}}}\|g\|^{2}_{L^{2}}+\nu  \epsilon\|\nabla g\|^{2}_{L^{2}}; \label{iii3.3}
\end{align}
\begin{align}
	b\left\langle|w|^{2\beta}w - |u|^{2\beta}u,g\right\rangle_{Y_{\beta}',Y_{\beta}} &\geq \frac{b}{2}\int_{\Omega}|g|^{2}(|w|^{2\beta}+|u|^{2\beta})dx \geq \frac{b}{2^{2\beta}}\||g|(|w|+|u|)^{\beta}\|^{2}_{L^{2}}. \label{iii3.1} 
\end{align}
By Mean Value Theorem, we have
\begin{align}
	\left|\left|\frac{l}{\nu}u\right|^{2\alpha} - \left|\frac{l}{\nu} u\right|^{2\beta}\right|\leq 2|\alpha-\beta|\left(\left|\frac{l}{\nu}u\right|^{2\alpha}+\left|\frac{l}{\nu}u\right|^{2\beta}\right)\left|\ln \left|\frac{l}{\nu}u\right|\right|. \nonumber
\end{align}


Henceforward, we divide in two cases for $\alpha$ and $\beta$:\\
\underline{Case $1<\alpha,\, \beta < 3$:}

\begin{align}
	\tilde{a}\frac{\nu}{l^{2}}\bigg\langle\left|\frac{l}{\nu}u\right|^{2\alpha}u-& \left|\frac{l}{\nu}u\right|^{2\beta}u,g\bigg\rangle_{Y_{\beta}',Y_{\beta}} \leq 2 \tilde{a}\frac{\nu^{2}}{l^{3}}|\alpha-\beta|\int_{\Omega}|g|\left(\left|\frac{l}{\nu}u\right|^{2\alpha+1}+\left|\frac{l}{\nu}u\right|^{2\beta+1}\right)\left|\ln \left|\frac{l}{\nu}u\right|\right| dx  \nonumber\\
	& \leq 4 \tilde{a}\frac{\nu^{2}}{l^{3}}|\alpha-\beta|\int_{\Omega}|g|\left(\frac{1}{3-\max\{\alpha,\,\beta\}}\left|\frac{l}{\nu}u\right|^{7}+\frac{1}{e}\left|\frac{l}{\nu}u\right|\right)dx 
	\nonumber\\
	& \leq \frac{4\tilde{a}|\alpha-\beta|}{3-\max\{\alpha,\,\beta\}}\frac{l^{4}}{\nu^{5}}\|g\|_{L^{6}}\||u|^{7}\|_{L^{\frac{6}{5}}}+\frac{4\tilde{a}|\alpha-\beta|}{e}\frac{\nu}{l^{2}}\|g\|_{L^{2}}\|u\|_{L^{2}}
	\nonumber\\
	& \leq \frac{4\tilde{a}C_{6}|\alpha-\beta|}{3-\max\{\alpha,\,\beta\}}\frac{l^{4}}{\nu^{5}}
	\left(\|\nabla g\|_{L^{2}}+\frac{1}{l}\|g\|_{L^{2}}\right)\|u\|^{7}_{L^\frac{42}{5}}+\frac{4\tilde{a}|\alpha-\beta|}{e}\frac{\nu}{l^{2}}\|g\|_{L^{2}}\|u\|_{L^{2}}\nonumber\\
	& \leq \frac{64\tilde{a}^{2}C^{2}_{6}|\alpha-\beta|^{2}}{(3-\max\{\alpha,\,\beta\})^{2}}\frac{l^{8}}{\nu^{10}}\left(\frac{1}{\nu}+\frac{1}{\tilde{\epsilon}\eta l^{2}}\right)C^{14}_{\frac{42}{5}}\left(\|\nabla u\|^{\frac{6}{7}}_{L^{2}}\|A u\|^{\frac{1}{7}}_{L^{2}}+\frac{1}{l^{\frac{8}{7}}}\|u\|_{L^{2}}\right)^{14}\nonumber\\
	&+
	\frac{32\tilde{a}^{2}|\alpha-\beta|^{2}}{\tilde{\epsilon}\eta e^{2}}\frac{\nu^{2}}{l^{4}}\|u\|^{2}_{L^{2}} +\frac{\nu}{4}\|\nabla g\|^{2}_{L^{2}} + \frac{\tilde{\epsilon}\eta}{4}\|g\|^{2}_{L^{2}}
	\nonumber\\
	& \leq \frac{2^{19}\tilde{a}^{2}C^{2}_{6}|\alpha-\beta|^{2}}{(3-\max\{\alpha,\,\beta\})^{2}}\frac{l^{8}}{\nu^{10}}\left(\frac{1}{\nu}+\frac{1}{\tilde{\epsilon}\eta l^{2}}\right)C^{14}_{\frac{42}{5}}\left(\|\nabla u\|^{12}_{L^{2}}\|A u\|^{2}_{L^{2}}+\frac{1}{l^{16}}\|u\|^{14}_{L^{2}}\right)\nonumber\\
	&+
	\frac{32\tilde{a}^{2}|\alpha-\beta|^{2}}{\tilde{\epsilon}\eta e^{2}}\frac{\nu^{2}}{l^{4}}\|u\|^{2}_{L^{2}} +\frac{\nu}{4}\|\nabla g\|^{2}_{L^{2}} + \frac{\tilde{\epsilon}\eta}{4}\|g\|^{2}_{L^{2}}.
	\label{iii3.4.0} 
\end{align}
Finally,
\begin{align}
	\left|\tilde{a}-\tilde{b}\right|\frac{\nu}{l^{2}}\left\langle \left|\frac{l}{\nu}u\right|^{2\beta}u,g\right\rangle_{Y_{\beta}',Y_{\beta}} &\leq \left|\tilde{a}-\tilde{b}\right|\frac{l^{4}}{\nu^{5}}\|g\|_{L^{6}}\left\||u|^{7}\right\|_{L^{\frac{6}{5}}} + \left|\tilde{a}-\tilde{b}\right|\frac{\nu}{l^{2}}\|g\|_{L^{2}}\left\|u\right\|_{L^{2}} \nonumber \\&\leq
	\left|\tilde{a}-\tilde{b}\right|\frac{l^{4}}{\nu^{5}}\left(\|\nabla g\|_{L^{2}}+\frac{1}{l}\|g\|_{L^{2}}\right)\left\|u\right\|_{L^{\frac{42}{5}}}^{7} + \left|\tilde{a}-\tilde{b}\right|\frac{\nu}{l^{2}}\|g\|_{L^{2}}\left\|u\right\|_{L^{2}}  \nonumber\\&\leq
	\left|\tilde{a}-\tilde{b}\right|^{2}\frac{l^{8}}{\nu^{10}}\left(\frac{1}{\nu}+\frac{2}{\bar{\epsilon}\eta l^{2}}\right)C^{14}_{\frac{42}{5}}\left(\|\nabla u\|^{\frac{6}{7}}_{L^{2}}\|A u\|^{\frac{1}{7}}_{L^{2}}+\frac{1}{l^{\frac{8}{7}}}\|u\|_{L^{2}}\right)^{14} \nonumber\\&+ \frac{2}{\bar{\epsilon}\eta}\left|\tilde{a}-\tilde{b}\right|^{2}\frac{\nu^{2}}{l^{4}}\|u\|^{2}_{L^{2}} +\frac{\nu}{4}\|\nabla g\|^{2}_{L^{2}}+ \frac{\bar{\epsilon}\eta}{4}\|g\|^{2}_{L^{2}} 
	\nonumber\\&\leq
	2^{14}\left|\tilde{a}-\tilde{b}\right|^{2}\frac{l^{8}}{\nu^{10}}\left(\frac{1}{\nu}+\frac{2}{\bar{\epsilon}\eta l^{2}}\right)C^{14}_{\frac{42}{5}}\left(\|\nabla u\|^{12}_{L^{2}}\|A u\|^{2}_{L^{2}}+\frac{1}{l^{16}}\|u\|^{14}_{L^{2}}\right) \nonumber\\&+ \frac{2}{\bar{\epsilon}\eta}\left|\tilde{a}-\tilde{b}\right|^{2}\frac{\nu^{2}}{l^{4}}\|u\|^{2}_{L^{2}} +\frac{\nu}{4}\|\nabla g\|^{2}_{L^{2}}+ \frac{\bar{\epsilon}\eta}{4}\|g\|^{2}_{L^{2}}. \label{iii3.4.0.1}
\end{align}

From (\ref{in11}) and the above inequalities, and choosing  $\widehat{\epsilon}= \tilde{\epsilon}=\bar{\epsilon}=1$ we obtain
\begin{align}
	\frac{1}{2}\frac{d}{dt} \|g\|_{L^{2}}^{2}&+\left[\nu\left(\frac{1}{2}-\epsilon\right)- \eta c_{o}h^{2}\right]\|\nabla g\|_{L^{2}}^{2}+ \left[\frac{b}{2^{2\beta}}-  \frac{1}{\nu^{\beta}(\eta\check{\epsilon})^{\beta-1}}\right]\||g|(|w|+|u|)^{\beta}\|^{2}_{L^{2}} \nonumber\\&\leq \eta\left[ \frac{(\beta-1)\check{\epsilon}}{\beta(16 \epsilon)^{\frac{\beta}{\beta-1}}} -\frac{1}{4}\right]\|g\|^{2}_{L^{2}} 
	\nonumber\\&+
	\frac{2^{19}\tilde{a}^{2}C^{2}_{6}|\alpha-\beta|^{2}}{(3-\max\{\alpha,\,\beta\})^{2}}\frac{l^{8}}{\nu^{10}}\left(\frac{1}{\nu}+\frac{1}{\eta l^{2}}\right)C^{14}_{\frac{42}{5}}\left(\|\nabla u\|^{12}_{L^{2}}\|A u\|^{2}_{L^{2}}+\frac{1}{l^{16}}\|u\|^{14}_{L^{2}}\right)+
	\frac{32\tilde{a}^{2}|\alpha-\beta|^{2}}{\eta e^{2}}\frac{\nu^{2}}{l^{4}}\|u\|^{2}_{L^{2}}
	\nonumber\\&+
	2^{14}\left|\tilde{a}-\tilde{b}\right|^{2}\frac{l^{8}}{\nu^{10}}\left(\frac{1}{\nu}+\frac{2}{\eta l^{2}}\right)C^{14}_{\frac{42}{5}}\left(\|\nabla u\|^{12}_{L^{2}}\|A u\|^{2}_{L^{2}}+\frac{1}{l^{16}}\|u\|^{14}_{L^{2}}\right) + \frac{2}{\eta}\left|\tilde{a}-\tilde{b}\right|^{2}\frac{\nu^{2}}{l^{4}}\|u\|^{2}_{L^{2}}. \nonumber
\end{align}%
Choosing, $\epsilon=\frac{1}{4}$, $\check{\epsilon}=\frac{\beta 4^{\frac{\beta}{\beta-1}}}{8(\beta-1)}$,  $\eta$  and $h$ as in the statement of Theorem \ref{CDA3a1.b}, $M_{1}$ and $M_{2}$ as given in Corollary \ref{da0101.4565.1}, we have
\begin{align}
	\frac{1}{2}\frac{d}{dt} \|g\|_{L^{2}}^{2}&+ \frac{\eta}{8}\|g\|^{2}_{L^{2}} \leq 
	\frac{2^{19}\tilde{a}^{2}C^{2}_{6}|\alpha-\beta|^{2}}{(3-\max\{\alpha,\,\beta\})^{2}}\frac{l^{8}}{\nu^{10}}\left(\frac{1}{\nu}+\frac{1}{\eta l^{2}}\right)C^{14}_{\frac{42}{5}}\left(M_{2}^{6}\|A u\|^{2}_{L^{2}}+\frac{1}{l^{16}}M_{1}^{7}\right)\nonumber\\& +
	\frac{32\tilde{a}^{2}|\alpha-\beta|^{2}}{\eta e^{2}}\frac{\nu^{2}}{l^{4}}M_{1}
	+
	2^{14}\left|\tilde{a}-\tilde{b}\right|^{2}\frac{l^{8}}{\nu^{10}}\left(\frac{1}{\nu}+\frac{2}{\eta l^{2}}\right)C^{14}_{\frac{42}{5}}\left(M_{2}^{6}\|A u\|^{2}_{L^{2}}+\frac{1}{l^{16}}M_{1}^{7}\right) \nonumber\\& + \frac{2}{\eta}\left|\tilde{a}-\tilde{b}\right|^{2}\frac{\nu^{2}}{l^{4}}M_{1}, \label{ngngnhgntr56} 
\end{align}%
for all $t\geq 0$. Using Gronwall's inequality, we conclude that
\begin{align}
	\|g(t)\|_{L^{2}}^{2}&\leq e^{-\frac{\eta}{8}t}\|g(0)\|_{L^{2}}^{2}+
	\frac{2^{20}\tilde{a}^{2}C^{2}_{6}|\alpha-\beta|^{2}}{(3-\max\{\alpha,\,\beta\})^{2}}\frac{l^{8}}{\nu^{10}}\left(\frac{1}{\nu}+\frac{1}{\eta l^{2}}\right)C^{14}_{\frac{42}{5}}\left(M_{2}^{6}\int_{0}^{t}e^{\frac{\eta}{8}(s-t)}\|A u(s)\|^{2}_{L^{2}}ds+\frac{8}{\eta l^{16}}M_{1}^{7}\right)\nonumber\\& +
	2^{15}\left|\tilde{a}-\tilde{b}\right|^{2}\frac{l^{8}}{\nu^{10}}\left(\frac{1}{\nu}+\frac{2}{\eta l^{2}}\right)C^{14}_{\frac{42}{5}}\left(M_{2}^{6}\int_{0}^{t}e^{\frac{\eta}{8}(s-t)}\|A u(s)\|^{2}_{L^{2}}ds+\frac{8}{\eta l^{16}}M_{1}^{7}\right) \nonumber\\& +\frac{2^{9}\tilde{a}^{2}|\alpha-\beta|^{2}}{\eta^{2} e^{2}}\frac{\nu^{2}}{l^{4}}M_{1}
	+ \frac{2^{5}}{\eta^{2}}\left|\tilde{a}-\tilde{b}\right|^{2}\frac{\nu^{2}}{l^{4}}M_{1}
	\nonumber \\ &\leq
	e^{-\frac{\eta}{8}t}\|g(0)\|_{L^{2}}^{2} +\frac{2^{9}\tilde{a}^{2}|\alpha-\beta|^{2}}{\eta^{2} e^{2}}\frac{\nu^{2}}{l^{4}}M_{1}
	+ \frac{2^{5}}{\eta^{2}}\left|\tilde{a}-\tilde{b}\right|^{2}\frac{\nu^{2}}{l^{4}}M_{1} \nonumber \\ &+
	\frac{2^{20}\tilde{a}^{2}C^{2}_{6}|\alpha-\beta|^{2}}{(3-\max\{\alpha,\,\beta\})^{2}}\frac{l^{8}}{\nu^{10}}\left(\frac{1}{\nu}+\frac{1}{\eta l^{2}}\right)C^{14}_{\frac{42}{5}}\left[M_{2}^{6}\left(\frac{4}{\nu}M_{2}+ \frac{16}{\eta \nu^{\frac{3\alpha-2}{\alpha-1}} a^{\frac{1}{\alpha-1}}} M_{2} +\frac{32}{\eta \nu^{2}}\|f\|^{2}_{L^{\infty}_{t}L^{2}}\right)+\frac{8}{\eta l^{16}}M_{1}^{7}\right]\nonumber\\& +
	2^{15}\left|\tilde{a}-\tilde{b}\right|^{2}\frac{l^{8}}{\nu^{10}}\left(\frac{1}{\nu}+\frac{2}{\eta l^{2}}\right)C^{14}_{\frac{42}{5}}\left[M_{2}^{6}\left(\frac{4}{\nu}M_{2}+ \frac{16}{\eta \nu^{\frac{3\alpha-2}{\alpha-1}} a^{\frac{1}{\alpha-1}}} M_{2} +\frac{32}{\eta \nu^{2}}\|f\|^{2}_{L^{\infty}_{t}L^{2}}\right)+\frac{8}{\eta l^{16}}M_{1}^{7}\right]. \nonumber
\end{align}%
Thus, the desired inequality in $H$-norm stated in Theorem \ref{CDA3a1.b} for the case $2\leq \alpha < 3$ or $2\leq \beta < 3$ is obtained.


\underline{Case  $1<\alpha,\, \beta < 2$:}


We choose $\widehat{\epsilon}=1$, $\epsilon =\frac{1}{4}$ and $\check{\epsilon}=\frac{\beta 4^{\frac{\beta}{\beta-1}}}{ 8(\beta-1)}$ in inequalities (\ref{iii3})-(\ref{iii3.3}). Since $1<\alpha, \beta <2$, we replace inequalities (\ref{iii3.4.0}) and (\ref{iii3.4.0.1}) by the following ones:
\begin{align}
	\tilde{a}\frac{\nu}{l^{2}}\left\langle\left|\frac{l}{\nu}u\right|^{2\alpha}u- \left|\frac{l}{\nu}u\right|^{2\beta}u,g\right\rangle_{Y_{\beta}',Y_{\beta}} &\leq 2 \tilde{a}\frac{\nu^{2}}{l^{3}}|\alpha-\beta|\int_{\Omega}|g|\left(\left|\frac{l}{\nu}u\right|^{2\alpha+1}+\left|\frac{l}{\nu}u\right|^{2\beta+1}\right)\left|\ln \left|\frac{l}{\nu}u\right|\right| dx  \nonumber\\& \leq 4 \tilde{a}\frac{\nu^{2}}{l^{3}}|\alpha-\beta|\int_{\Omega}|g|\left(\frac{1}{2-\max\{\alpha,\,\beta\}}\left|\frac{l}{\nu}u\right|^{5}+\frac{1}{e}\left|\frac{l}{\nu}u\right|\right)dx 
	\nonumber\\& \leq \frac{4\tilde{a}|\alpha-\beta|}{2-\max\{\alpha,\,\beta\}}\frac{l^{2}}{\nu^{3}}\|g\|_{L^{6}}\||u|^{5}\|_{L^{\frac{6}{5}}}+\frac{4\tilde{a}|\alpha-\beta|}{e}\frac{\nu}{l^{2}}\|g\|_{L^{2}}\|u\|_{L^{2}}
	\nonumber\\& \leq \frac{4\tilde{a}C_{6}|\alpha-\beta|}{2-\max\{\alpha,\,\beta\}}\frac{l^{2}}{\nu^{3}}\left(\|\nabla g\|_{L^{2}}+\frac{1}{l}\|g\|_{L^{2}}\right)\|u\|^{5}_{L^{6}}\nonumber\\&+\frac{4\tilde{a}|\alpha-\beta|}{e}\frac{\nu}{l^{2}}\|g\|_{L^{2}}\|u\|_{L^{2}} 
	\nonumber\\& \leq \frac{64\tilde{a}^{2}C^{2}_{6}|\alpha-\beta|^{2}}{(2-\max\{\alpha,\,\beta\})^{2}}\frac{l^{4}}{\nu^{6}}\left(\frac{1}{\nu}+\frac{1}{\tilde{\epsilon}\eta l^{2}}\right)C^{10}_{6}\left(\|\nabla u\|^{2}_{L^{2}}+\frac{1}{l^{2}}\|u\|^{2}_{L^{2}}\right)^{5}\nonumber\\&+
	\frac{32\tilde{a}^{2}|\alpha-\beta|^{2}}{\tilde{\epsilon}\eta e^{2}}\frac{\nu^{2}}{l^{4}}\|u\|^{2}_{L^{2}} +\frac{\nu}{4}\|\nabla g\|^{2}_{L^{2}} + \frac{\tilde{\epsilon}\eta}{4}\|g\|^{2}_{L^{2}};
	\label{iii3.4.0.a.a} 
\end{align}
and
\begin{align}
	\left|\tilde{a}-\tilde{b}\right|\frac{\nu}{l^{2}}\left\langle \left|\frac{l}{\nu}u\right|^{2\beta}u,g\right\rangle_{Y_{\beta}',Y_{\beta}} &\leq \left|\tilde{a}-\tilde{b}\right|\frac{l^{2}}{\nu^{3}}\|g\|_{L^{6}}\left\||u|^{5}\right\|_{L^{\frac{6}{5}}} + \left|\tilde{a}-\tilde{b}\right|\frac{\nu}{l^{2}}\|g\|_{L^{2}}\left\|u\right\|_{L^{2}} \nonumber \\&\leq
	\left|\tilde{a}-\tilde{b}\right|\frac{l^{2}}{\nu^{3}}\left(\|\nabla g\|_{L^{2}}+\frac{1}{l}\|g\|_{L^{2}}\right)\left\|u\right\|_{L^{6}}^{5} + \left|\tilde{a}-\tilde{b}\right|\frac{\nu}{l^{2}}\|g\|_{L^{2}}\left\|u\right\|_{L^{2}}  \nonumber\\&\leq
	\left|\tilde{a}-\tilde{b}\right|^{2}\frac{l^{4}}{\nu^{6}}\left(\frac{1}{\nu}+\frac{2}{\bar{\epsilon}\eta l^{2}}\right)\|u\|_{L^{6}}^{10} + \frac{2}{\bar{\epsilon}\eta}\left|\tilde{a}-\tilde{b}\right|^{2}\frac{\nu^{2}}{l^{4}}\|u\|^{2}_{L^{2}} +\frac{\nu}{4}\|\nabla g\|^{2}_{L^{2}}+ \frac{\bar{\epsilon}\eta}{4}\|g\|^{2}_{L^{2}} 
	\nonumber\\&\leq
	2\left|\tilde{a}-\tilde{b}\right|^{2}\frac{l^{4}}{\nu^{6}}\left(\frac{1}{\nu}+\frac{2}{\bar{\epsilon}\eta l^{2}}\right)C^{10}_{6}\left(\|\nabla u\|^{2}_{L^{2}}+\frac{1}{l^{2}}\|u\|^{2}_{L^{2}}\right)^{5} + \frac{2}{\bar{\epsilon}\eta}\left|\tilde{a}-\tilde{b}\right|^{2}\frac{\nu^{2}}{l^{4}}\|u\|^{2}_{L^{2}} \nonumber\\&+\frac{\nu}{4}\|\nabla g\|^{2}_{L^{2}}+ \frac{\bar{\epsilon}\eta}{4}\|g\|^{2}_{L^{2}}. \label{iii3.4.0.1.a.a}
\end{align}
Then, choosing $\tilde{\epsilon}=\bar{\epsilon}=1$, using (\ref{in11}), hypothesis (\ref{12ee2}),  $M_{1}$ and $M_{2}$ as given in Corollary \ref{da0101.4565.1}, we obtain
\begin{align}
	\frac{1}{2}\frac{d}{dt} \|g\|_{L^{2}}^{2}&+ \frac{\eta}{8}\|g\|^{2}_{L^{2}} \leq \frac{64\tilde{a}^{2}C^{12}_{6} |\alpha-\beta|^{2}}{(2-\max\{\alpha,\,\beta\})^{2}}\frac{l^{2}(\eta l^{2}+2\nu)}{\eta\nu^{7}}\left(M_{2}+\frac{1}{l^{2}}M_{1}\right)^{5}\nonumber\\&+
	\frac{32\tilde{a}^{2}|\alpha-\beta|^{2}}{e^{2}}\frac{\nu^{2}}{\eta l^{4}}M_{1}   
	\nonumber\\&+
	2C^{10}_{6}\left|\tilde{a}-\tilde{b}\right|^{2}\frac{l^{2}(\eta l^{2}+2\nu)}{\eta\nu^{7}} \left(M_{2}+\frac{1}{l^{2}}M_{1}\right)^{5} + 2\left|\tilde{a}-\tilde{b}\right|^{2}\frac{\nu^{2}}{\eta l^{4}}M_{1},\,\, \forall\, t \geq 0.
\end{align}%
By Gronwall's inequality,  we conclude the inequality stated in Theorem \ref{CDA3a1.b} in $H$-norm,  for the case $1< \alpha, \beta < 2$. 

\begin{flushright}
	\rule{2mm}{2mm}
\end{flushright}


\section{Proof of Theorem \ref{CDA3a1.b.c.d}} \label{sec6}

We add a convenient term in inequality (\ref{iii3}) and get
\begin{align}
	\left|\eta (I_{h}(g)-g,g)_{L^{2}}\right| \leq \eta\widehat{\epsilon}\|I_{h}(g)-g\|^{2}_{L^{2}} +  \eta\frac{1}{4\widehat{\epsilon}}\|g\|^{2}_{L^{2}}\leq  \eta c_{o}h^{2}\widehat{\epsilon}\|\nabla g\|_{L^{2}}^{2} +\eta c_{1}h^{4}\widehat{\epsilon}\|A g\|_{L^{2}}^{2}+ \eta\frac{1}{4\widehat{\epsilon}}\|g\|^{2}_{L^{2}}. \label{iii3ull0}
\end{align}
Choose now $\epsilon =\frac{1}{4}$ and $\check{\epsilon}=\frac{\beta 4^{\frac{\beta}{\beta-1}}}{ 8(\beta-1)}$ in inequalities (\ref{iii3})-(\ref{iii3.3}),  $\tilde{\epsilon}=\bar{\epsilon}=\frac{\nu}{\eta l^{2}}$ in (\ref{iii3.4.0.a.a})-(\ref{iii3.4.0.1.a.a}) and $\widehat{\epsilon}=1$ in (\ref{iii3ull0}). By (\ref{in11}), hypotheses (\ref{12ee2}) and division by $l^{2}$, we obtain
\begin{align}
	\frac{1}{2}\frac{d}{dt} \frac{1}{l^{2}}\|g\|_{L^{2}}^{2}  &\leq \frac{1}{l^{2}}\|g\|^{2}_{L^{2}}\left[- \frac{3\eta}{8}+ \frac{\nu}{2l^{2}}\right]+ \frac{\eta c_{1}h^{4}}{l^{2}}\|A g\|_{L^{2}}^{2} +\nonumber\\&+
	|\alpha-\beta|^{2}\left[\frac{128\tilde{a}^{2}C^{12}_{6}}{(2-\max\{\alpha,\,\beta\})^{2}}\frac{l^{2}}{\nu^{7}}\left(\|\nabla u\|^{2}_{L^{2}}+\frac{1}{l^{2}}\|u\|^{2}_{L^{2}}\right)^{5}\nonumber+\frac{32\tilde{a}^{2}}{ e^{2}}\frac{\nu}{l^{4}}\|u\|^{2}_{L^{2}}\right] 
	\nonumber\\&+
	\left|\tilde{a}-\tilde{b}\right|^{2}\left[6C^{10}_{6}\frac{l^{2}}{\nu^{7}}\left(\|\nabla u\|^{2}_{L^{2}}+\frac{1}{l^{2}}\|u\|^{2}_{L^{2}}\right)^{5} + 2\frac{\nu}{l^{4}}\|u\|^{2}_{L^{2}}\right] + \frac{\eta c_{1}h^{4}}{l^{2}}\|A g\|_{L^{2}}^{2}. \label{ngngnhgntr56uiyi}
\end{align}

Multiplying the system (\ref{eq111.0}) by $Ag$, integrating over $\Omega$, we have
\begin{align}\label{in1123.1}
	\frac{1}{2}\frac{d}{dt} \|\nabla g\|_{L^{2}}^{2}&+\nu\|A g\|_{L^{2}}^{2}+ b\left(|w|^{2\beta}w - |u|^{2\beta}u,Ag\right)_{L^{2}}  = 
	\tilde{a}\frac{\nu}{l^{2}}\left(\left|\frac{l}{\nu}u\right|^{2\alpha}u-\left|\frac{l}{\nu}u\right|^{2\beta}u, Ag\right)_{L_{2}}\nonumber \\&+
	\left(\tilde{a}-\tilde{b}\right)\frac{\nu}{l^{2}}\left(\left|\frac{l}{\nu}u\right|^{2\beta}u, Ag\right)_{L_{2}}
	-\left(g\cdot \nabla u, A g\right)_{L^{2}}\\&- \left(w\cdot \nabla g, A g\right)_{L^{2}}-\eta\left(I_{h}(g)-g,Ag\right)_{L^{2}}-\eta\|\nabla g\|^{2}_{L^{2}}.\nonumber	
\end{align}%
In order to estimate the terms of (\ref{in1123.1}), we use (\ref{0in1}), (\ref{gn}), (\ref{dam02}) and Hölder and Young's inequalities:
\begin{align}\label{in1123.2}
	|\eta(I_{h}(g)-g,Ag)_{L^{2}}| &\leq \eta \left(\sqrt{c_{0}h^{2}\|\nabla g\|^{2}_{L^{2}}+ c_{1}h^{4}\|A g\|^{2}_{L^{2}}}\right)\|Ag\|_{L^{2}}\nonumber\\
	&  \leq \frac{8\eta^{2}c_{0}h^{2}}{\nu}\|\nabla g\|^{2}_{L^{2}}+ \frac{8\eta^{2}c_{1}h^{4}}{\nu}\|Ag\|^{2}_{L^{2}}+ \frac{\nu}{32}\|Ag\|^{2}_{L^{2}};
\end{align}
\begin{align}\label{xx}
	\left|(g\cdot \nabla u, A g)_{L^{2}}\right| &\leq \|g\|_{L^{\infty}}\|\nabla u\|_{L^{2}}\|Ag\|_{L^{2}} \nonumber\\&\leq
	C_{\infty}\left(\|\nabla g\|^{\frac{1}{2}}_{L^{2}}\|Ag\|^{\frac{1}{2}}_{L^{2}}+\frac{1}{l^{\frac{3}{2}}}\|g\|_{L^{2}}\right)\|\nabla u\|_{L^{2}}\|Ag\|_{L^{2}}
	\nonumber\\&= C_{\infty}\|\nabla g\|^{\frac{1}{2}}_{L^{2}}\|Ag\|^{\frac{3}{2}}_{L^{2}}\|\nabla u\|_{L^{2}}+
	C_{\infty}\frac{1}{l^{\frac{3}{2}}}\|g\|_{L^{2}}\|\nabla u\|_{L^{2}}\|Ag\|_{L^{2}}
	\nonumber\\&\leq  \frac{432}{\nu^{3}}C^{4}_{\infty}\|\nabla g\|^{2}_{L^{2}}\|\nabla u\|^{4}_{L^{2}}+
	\frac{4}{\nu}C^{2}_{\infty}\frac{1}{l^{3}}\|g\|^{2}_{L^{2}}\|\nabla u\|^{2}_{L^{2}};
\end{align}
\begin{align}\label{in1123.3}
	\left|(w\cdot \nabla g, A g)_{L^{2}}\right| &\leq \|w\|_{L^{6}}\|\nabla g\|_{L^{3}}\|A g\|_{L^{2}} 
	\nonumber\\&\leq C_{6}C_{3}\left(\| \nabla  w\|_{L^{2}}+ \frac{1}{l}\| w\|_{L^{2}}\right)\left(\|\nabla g\|^{\frac{1}{2}}_{L^{2}}\| A g\|^{\frac{1}{2}}_{L^{2}}+\frac{1}{l^{\frac{1}{2}}}\|\nabla g\|_{L^{2}}\right)\|A g\|_{L^{2}}
	\nonumber\\&= C_{6}C_{3}\left(\| \nabla  w\|_{L^{2}}+ \frac{1}{l}\| w\|_{L^{2}}\right)\frac{1}{l^{\frac{1}{2}}}\|\nabla g\|_{L^{2}}\|A g\|_{L^{2}} \\&+
	C_{6}C_{3}\left(\| \nabla  w\|_{L^{2}}+ \frac{1}{l}\| w\|_{L^{2}}\right)\|\nabla g\|^{\frac{1}{2}}_{L^{2}}\|A g\|^{\frac{3}{2}}_{L^{2}}
	\nonumber \\&\leq
	\frac{8}{\nu}C^{2}_{6}C^{2}_{3}\left(\| \nabla  w\|^{2}_{L^{2}}+ \frac{1}{l^{2}}\| w\|^{2}_{L^{2}}\right)\frac{1}{l}\|\nabla g\|^{2}_{L^{2}} +
	\frac{12^{3}}{\nu^{3}} C^{4}_{6}C^{4}_{3}\left(\| \nabla  w\|^{2}_{L^{2}}+ \frac{1}{l^{2}}\| w\|^{2}_{L^{2}}\right)^{2}\|\nabla g\|^{2}_{L^{2}}
	\nonumber \\&+\frac{\nu}{8}\|A g\|^{2}_{L^{2}}. \nonumber
\end{align}

Since $1 <  \beta  < 2$, denoting $\displaystyle J_{10}= b\left|(|w|^{2\beta}w - |u|^{2\beta}u,Ag)_{L^{2}}\right|$ and using (\ref{gn}) and (\ref{dam02}),  we get

\begin{align}\label{in1123.7}
	J_{10}&\leq |b|\kappa(2\beta)\int_{\Omega}|g||Ag|(|u|+|v|)^{2\beta}dx \leq |b|\kappa(2\beta)\|g\|_{L^{6}}\|Ag\|_{L^{2}}\|(u+w)^{2\beta}\|_{L^{3}} \nonumber\\
	&\leq   8|b|\kappa(2\beta)C_{6}C^{2\beta}_{6\beta}\left(\|\nabla g\|_{L^{2}}+ \frac{1}{l}\|g\|_{L^{2}}\right)\|Ag\|_{L^{2}}\left(\|\nabla (u+w) \|^{1+\beta}_{L^{2}}\|A (u+w)\|^{\beta-1}_{L^{2}}+ \frac{1}{l^{3\beta-1}}\|u+w\|^{2\beta}_{L^{2}}\right)
	\nonumber\\
	&\leq   8|b|\kappa(2\beta)C_{6}C^{2\beta}_{6\beta}\left(\|\nabla g\|_{L^{2}}+ \frac{1}{l}\|g\|_{L^{2}}\right)\|Ag\|_{L^{2}}\nonumber\\
	&\cdot\left(\|\nabla (u+w) \|_{L^{2}}^{1+\beta}\left(\|A g\|_{L^{2}} + 2\|A u\|_{L^{2}} \right)^{\beta-1}+ \frac{1}{l^{3\beta-1}}\|u+w\|^{2\beta}_{L^{2}}\right)
	\nonumber\\
	&\leq  16|b|\kappa(2\beta)C_{6}C^{2\beta}_{6\beta}\left(\|\nabla g\|_{L^{2}}+ \frac{1}{l}\|g\|_{L^{2}}\right)\|Ag\|_{L^{2}}\left(\frac{1}{l^{3\beta-1}}\|u+w\|^{2\beta}_{L^{2}}+\|\nabla (u+w) \|_{L^{2}}^{1+\beta}\|A u\|_{L^{2}}^{\beta-1}\right)\nonumber\\
	& + 8|b|\kappa(2\beta)C_{6}C^{2\beta}_{6\beta}\left(\|\nabla g\|_{L^{2}}+ \frac{1}{l}\|g\|_{L^{2}}\right)\|Ag\|^{\beta}_{L^{2}}\|\nabla (u+w)\|_{L^{2}}^{1+\beta}
	\nonumber\\
	&\leq   \frac{2^{12}}{\nu}|b|^{2}\kappa^{2}(2\beta)C^{2}_{6}C^{4\beta}_{6\beta}\left(\|\nabla g\|^{2}_{L^{2}}+ \frac{1}{l^{2}}\|g\|^{2}_{L^{2}}\right)\left(\frac{1}{l^{6\beta-2}}\|u+w\|^{4\beta}_{L^{2}}+\|A u\|_{L^{2}}^{2\beta-2}\|\nabla (u+w) \|_{L^{2}}^{2+2\beta}\right)
	\nonumber\\
	&+\left(\frac{16}{\nu}\right)^{\frac{\beta}{2-\beta}}2^{\frac{7}{2-\beta}}|b|^{\frac{2}{2-\beta}}\kappa^{\frac{2}{2-\beta}}(2\beta)C^{\frac{2}{2-\beta}}_{6}C^{\frac{4\beta}{2-\beta}}_{6\beta}\left(\|\nabla g\|^{2}_{L^{2}}+ \frac{1}{l^{2}}\|g\|^{2}_{L^{2}}\right)^{\frac{1}{2-\beta}}\|\nabla u+w \|_{L^{2}}^{\frac{2+2\beta}{2-\beta}}  +  \frac{\nu}{8}\|Ag\|^{2}_{L^{2}}
	\nonumber\\
	&\leq   \frac{2^{12}}{\nu}|b|^{2}\kappa^{2}(2\beta)C^{2}_{6}C^{4\beta}_{6\beta}\left(\|\nabla g\|^{2}_{L^{2}}+ \frac{1}{l^{2}}\|g\|^{2}_{L^{2}}\right)\left(\frac{1}{l^{6\beta-2}}\|u+w\|^{4\beta}_{L^{2}}+\frac{l}{b^{2}}\|A u\|_{L^{2}}^{2}+\frac{b^{\frac{2\beta-2}{2-\beta}}}{l^{\frac{\beta-1}{2-\beta}}}\|\nabla (u+w) \|_{L^{2}}^{\frac{2+2\beta}{2-\beta}}\right)
	\nonumber\\
	&+
	\left(\frac{16}{\nu}\right)^{\frac{\beta}{2-\beta}}2^{\frac{7}{2-\beta}}|b|^{\frac{2}{2-\beta}}\kappa^{\frac{2}{2-\beta}}(2\beta)C^{\frac{2}{2-\beta}}_{6}C^{\frac{4\beta}{2-\beta}}_{6\beta}\left(\|\nabla g\|^{2}_{L^{2}}+ \frac{1}{l^{2}}\|g\|^{2}_{L^{2}}\right)^{\frac{1}{2-\beta}}\|\nabla u+w \|_{L^{2}}^{\frac{2+2\beta}{2-\beta}}  +  \frac{\nu}{8}\|Ag\|^{2}_{L^{2}}.  
	\end{align}
Moreover, for $\displaystyle J_{11}= \tilde{a}\frac{\nu}{l^{2}}\left(\left|\frac{l}{\nu}u\right|^{2\alpha}u-\left|\frac{l}{\nu}u\right|^{2\beta}u, Ag\right)_{L_{2}}$, we obtain
\begin{align}\label{iii3.4}
	J_{11}&\leq 2\tilde{a}\frac{\nu^{2}}{l^{3}}|\alpha-\beta|\int_{\Omega}|Ag|\left(\left|\frac{l}{\nu}u\right|^{2\alpha+1}+\left|\frac{l}{\nu}u\right|^{2\beta+1}\right)\left|\ln \left|\frac{l}{\nu}u\right|\right| dx  \nonumber\\& \leq 4\tilde{a}\frac{\nu^{2}}{l^{3}}|\alpha-\beta|\int_{\Omega}|Ag|\left(\frac{1}{2-\max\{\alpha,\,\beta\}}\left|\frac{l}{\nu}u\right|^{5}+\frac{1}{e}\left|\frac{l}{\nu}u\right|\right)dx 
	\nonumber\\& \leq \frac{4\tilde{a}|\alpha-\beta|}{2-\max\{\alpha,\,\beta\}}\frac{l^{2}}{\nu^{3}}\|Ag\|_{L^{2}}\|u\|^{5}_{L^{10}}+\frac{4\tilde{a}|\alpha-\beta|}{e}\frac{\nu}{l^{2}}\|Ag\|_{L^{2}}\|u\|_{L^{2}}\\ 
		 & \leq \frac{32\tilde{a}^{2}|\alpha-\beta|^{2}}{(2-\max\{\alpha,\,\beta\})^{2}}\frac{l^{4}}{\nu^{7}}\|u\|^{10}_{L^{10}}+
	\frac{32\tilde{a}^{2}|\alpha-\beta|^{2}}{e^{2}} \frac{\nu}{l^{4}}\|u\|^{2}_{L^{2}} \ +\frac{\nu}{4}\|Ag\|^{2}_{L^{2}}
	\nonumber\\& \leq \frac{32\tilde{a}^{2}|\alpha-\beta|^{2}}{(2-\max\{\alpha,\,\beta\})^{2}}\frac{l^{4}}{\nu^{7}}C^{10}_{10}\left(\|\nabla u\|^{\frac{4}{5}}_{L^{2}} \|Au\|^{\frac{1}{5}}_{L^{2}}+\frac{1}{l^{\frac{6}{5}}}\|u\|_{L^{2}}\right)^{10}+
	\frac{32\tilde{a}^{2}|\alpha-\beta|^{2}}{e^{2}}\frac{\nu}{l^{4}}\|u\|^{2}_{L^{2}} +\frac{\nu}{4}\|Ag\|^{2}_{L^{2}}
	\nonumber\\& \leq \frac{2^{14}\tilde{a}^{2}|\alpha-\beta|^{2}}{(2-\max\{\alpha,\,\beta\})^{2}}C^{10}_{10}\frac{l^{4}}{\nu^{7}}\left(\|\nabla u\|^{8}_{L^{2}} \|Au\|^{2}_{L^{2}}+\frac{1}{l^{12}}\|u\|^{10}_{L^{2}}\right)+
	\frac{32\tilde{a}^{2}|\alpha-\beta|^{2}}{e^{2}} \frac{\nu}{l^{4}}\|u\|^{2}_{L^{2}} +\frac{\nu}{4}\|Ag\|^{2}_{L^{2}}.\nonumber
\end{align}
Finally, we have
\begin{align}\label{xx1}
	\left|\left(\tilde{a}-\tilde{b}\right)\frac{\nu}{l^{2}}\left(\left|\frac{l}{\nu}u\right|^{2\beta}u, Ag\right)_{L_{2}}\right| &\leq |\tilde{a}-\tilde{b}|\frac{l^{2}}{\nu^{3}} \|u\|^{5}_{L^{10}}\|Ag\|_{L^{2}}+ |\tilde{a}-\tilde{b}|\frac{\nu}{l^{2}} \|u\|_{L^{2}}\|Ag\|_{L^{2}}
	\nonumber\\&\leq 2^{11}|\tilde{a}-\tilde{b}|^{2}\frac{l^{4}}{\nu^{7}}C_{10}^{10} \left(\|\nabla u\|^{8}_{L^{2}} \|Au\|^{2}_{L^{2}}+\frac{1}{l^{12}}\|u\|^{10}_{L^{2}}\right)\nonumber\\&+ 4|\tilde{a}-\tilde{b}|^{2}\frac{\nu}{l^{4}} \|u\|^{2}_{L^{2}}+\frac{\nu}{8}\|Ag\|^{2}_{L^{2}}.
\end{align}

Considering $M_{1}$ and $M_{2}$ given in Corollary \ref{da0101.4565.1}, we define the following quantities:
\begin{eqnarray}
	Z_{1}(t)&=&432C^{4}_{\infty} \frac{1}{\nu^{3}}M_{2}^{2}+
	4C^{2}_{\infty} \frac{1}{l \nu}M_{2}+
	8C^{2}_{6}C^{2}_{3} \frac{1}{\nu l}\left(\| \nabla  w(t)\|^{2}_{L^{2}}+ \frac{1}{l^{2}}\| w(t)\|^{2}_{L^{2}}\right) \nonumber\\
	&+&
	12^{3} C^{4}_{6}C^{4}_{3}\frac{1}{\nu^{3}}\left(\| \nabla  w(t)\|^{2}_{L^{2}}+ \frac{1}{l^{2}}\| w(t)\|^{2}_{L^{2}}\right)^{2}
	+2^{4\beta+11}\kappa^{2}(2\beta)C^{2}_{6}C^{4\beta}_{6\beta}\frac{b^{2}}{\nu l^{6\beta-2}}\left(M_{1}^{2\beta} + \|w(t)\|^{4\beta}_{L^{2}}\right)
	\nonumber\\
	&+&
	2^{\frac{24-9\beta}{2-\beta}}\kappa^{2}(2\beta)C^{2}_{6}C^{4\beta}_{6\beta}
	\frac{b^{\frac{2}{2-\beta}}}{l^{\frac{\beta-1}{2-\beta}}\nu}
	\left(M_{2}^{\frac{1+\beta}{2-\beta}}+\|\nabla w(t) \|_{L^{2}}^{\frac{2+2\beta}{2-\beta}}\right)
	\label{Z1}\\
	&+&
	2^{\frac{7\beta+7}{2-\beta}}\kappa^{\frac{2}{2-\beta}}(2\beta)C^{\frac{2}{2-\beta}}_{6}C^{\frac{4\beta}{2-\beta}}_{6\beta}
	\frac{b^{\frac{2}{2-\beta}}}{\nu^{\frac{\beta}{2-\beta}}}
	\left(\|\nabla g(t)\|^{2}_{L^{2}}+ \frac{1}{l^{2}}\|g(t)\|^{2}_{L^{2}}\right)^{\frac{\beta-1}{2-\beta}}\left(M_{2}^{\frac{1+\beta}{2-\beta}} + \|\nabla w(t) \|_{L^{2}}^{\frac{2+2\beta}{2-\beta}}\right)
	\nonumber\\
	 &+& \frac{\nu}{2l^{2}};\nonumber
\end{eqnarray}
\begin{equation}\label{z2ul}
Z_{2}(t)= 2^{12}\kappa^{2}(2\beta)C^{2}_{6}C^{4\beta}_{6\beta}\frac{l}{\nu}\|A u(t)\|_{L^{2}}^{2}; 
\end{equation}
\begin{equation}\label{Z3}
Z_{3}= \frac{64\tilde{a}^{2}}{e^{2}}\frac{\nu}{l^{4}}M_{1} +
\frac{2^{15}\tilde{a}^{2}(C^{12}_{6}+C^{10}_{10})}{(2-\max\{\alpha,\,\beta\})^{2}}\frac{l^{2}}{\nu^{7}}\left(M_{2}+\frac{1}{l^{2}}M_{1}\right)^{5};
\end{equation}
\begin{equation}\label{Z4}
Z_{4}(t)= \frac{2^{14}\tilde{a}^{2}C^{10}_{10}}{(2-\max\{\alpha,\,\beta\})^{2}} \frac{l^{4}}{\nu^{7}}
M_{2}^{4} \|Au(t)\|^{2}_{L^{2}};
\end{equation}
\begin{equation}\label{Z5}
Z_{5}= 6\frac{\nu}{l^{4}} M_{1}+
2^{12}(C^{10}_{6}+ C_{10}^{10})\frac{l^{2}}{\nu^{7}}\left(M_{2}+\frac{1}{l^{2}}M_{1}\right)^{5};
\end{equation}
\begin{equation}\label{z6ul}
Z_{6}(t)= 2^{11}C_{10}^{10}\frac{l^{4}}{\nu^{7}}M_{2}^{4} \|Au(t)\|^{2}_{L^{2}}. 
\end{equation}

Then, taking into account (\ref{ngngnhgntr56uiyi})-\eqref{in1123.1}, the estimates \eqref{in1123.2}, \eqref{xx}, \eqref{in1123.3}, \eqref{in1123.7}, \eqref{iii3.4} and \eqref{xx1} yield  
\begin{align}
	\frac{1}{2}\frac{d}{dt} \left(\|\nabla g(t)\|_{L^{2}}^{2}+ \frac{1}{l^{2}}\|g(t)\|_{L^{2}}^{2}\right) &\leq \left(\|\nabla g(t)\|_{L^{2}}^{2}+ \frac{1}{l^{2}}\|g(t)\|_{L^{2}}^{2}\right)\left(-\frac{3\eta}{8} + \frac{8\eta^{2}c_{0}h^{2}}{\nu} +Z_{1}(t)+Z_{2}(t)\right)\nonumber\\&+ |\alpha-\beta|^{2} \left(Z_{3}+Z_{4}(t)\right)+ |\tilde{a}-\tilde{b}|^{2}\left(Z_{5}+Z_{6}(t)\right),
	\label{zzxczxcfds}
\end{align}
for all $t\geq 0$. Besides, note that 
\begin{align}
	\|\nabla g(0)\|_{L^{2}}^{2}+ \frac{1}{l^{2}} \|g(0)\|_{L^{2}}^{2}\leq 2\|\nabla u(0)\|_{L^{2}}^{2}+ \frac{2}{l^{2}} \|u(0)\|_{L^{2}}^{2}+ 2\|\nabla w(0)\|_{L^{2}}^{2}+ \frac{2}{l^{2}} \|w(0)\|_{L^{2}}^{2} \leq 4M. \nonumber
\end{align}
Next, consider 
\begin{align}\label{uioudrt5e5}
	H&= M_{2} +\frac{1}{l^{2}}M_{1} + e^{2^{13}\kappa^{2}(2\beta)C^{2}_{6}C^{4\beta}_{6\beta}\frac{l}{\nu^{2}}M_{2}}\Bigg\{
	4M\nonumber\\&+2|\alpha-\beta|^{2} \left[\frac{2^{3}}{\eta}
	\left[\frac{64\tilde{a}^{2}}{e^{2}}\frac{\nu}{l^{4}}M_{1} +
	\frac{2^{15}\tilde{a}^{2}(C^{12}_{6}+C^{10}_{10})}{(2-\max\{\alpha,\,\beta\})^{2}}\frac{l^{2}}{\nu^{7}}\left(M_{2}+\frac{1}{l^{2}}M_{1}\right)^{5}\right] \right.\nonumber\\&\left.+\frac{2^{14}\tilde{a}^{2}C^{10}_{10}}{(2-\max\{\alpha,\,\beta\})^{2}} \frac{l^{4}}{\nu^{7}}
	M_{2}^{4} \left(\frac{4}{\nu}M_{2}+ \frac{16}{\eta \nu^{\frac{3\alpha-2}{\alpha-1}} a^{\frac{1}{\alpha-1}}} M_{2} +\frac{32}{\eta \nu^{2}}\|f\|^{2}_{L^{\infty}_{t}L^{2}}\right)\right]
	\nonumber\\&+
	2|\tilde{a}-\tilde{b}|^{2} \left[\frac{2^{3}}{\eta}
	\left[6\frac{\nu}{l^{4}} M_{1}+
	2^{12}(C^{10}_{6}+ C_{10}^{10})\frac{l^{2}}{\nu^{7}}\left(M_{2}+\frac{1}{l^{2}}M_{1}\right)^{5}\right]
	\right.\\&\left.+2^{11}C_{10}^{10}\frac{l^{4}}{\nu^{7}}M_{2}^{4} \left(\frac{4}{\nu}M_{2}+ \frac{16}{\eta \nu^{\frac{3\alpha-2}{\alpha-1}} a^{\frac{1}{\alpha-1}}} M_{2} +\frac{32}{\eta \nu^{2}}\|f\|^{2}_{L^{\infty}_{t}L^{2}}\right)\right]\Bigg\}.\nonumber 
\end{align}
and
\begin{align}
	T^{*}=\sup\left\{t^{*}>0;\,\, \|\nabla g(t)\|_{L^{2}}^{2}+ \frac{1}{l^{2}} \|g(t)\|_{L^{2}}^{2}\leq H,\,\, \forall\, 0\leq t\leq t^{*} \right\}. \label{ttt2}
\end{align}
Then,
\begin{align}
	\|\nabla w(t)\|_{L^{2}}^{2}+ \frac{1}{l^{2}} \|w(t)\|_{L^{2}}^{2}&\leq 2\|\nabla g(t)\|_{L^{2}}^{2}+ \frac{2}{l^{2}} \|g(t)\|_{L^{2}}^{2} + 2\|\nabla u(t)\|_{L^{2}}^{2}+ \frac{2}{l^{2}} \|u(t)\|_{L^{2}}^{2} \nonumber \\&\leq 2H+2M_{2} +\frac{2}{l^{2}}M_{1} \leq 4H,\nonumber
\end{align}
for all $0\leq t< T^{*}$. Defining
\begin{align} \label{dfbdfr54}
	\tilde{Z}_{1}&=432C^{4}_{\infty} \frac{1}{\nu^{3}}M_{2}^{2}+
	4C^{2}_{\infty} \frac{1}{\nu}M_{2}+
	8C^{2}_{6}C^{2}_{3} \frac{1}{\nu l}4H \nonumber\\&+
	12^{3} C^{4}_{6}C^{4}_{3}\frac{1}{\nu^{3}}16H^{2}
	+2^{4\beta+11}\kappa^{2}(2\beta)C^{2}_{6}C^{4\beta}_{6\beta}\frac{b^{2}}{\nu l^{6\beta-2}}\left(M_{1}^{2\beta} + (l^{2}4H)^{2\beta}_{L^{2}}\right)
	\nonumber\\&+
	2^{\frac{24-9\beta}{2-\beta}}\kappa^{2}(2\beta)C^{2}_{6}C^{4\beta}_{6\beta}
	\frac{b^{\frac{2}{2-\beta}}}{l^{\frac{\beta-1}{2-\beta}}\nu}
	\left(M_{2}^{\frac{1+\beta}{2-\beta}}+(4H)^{\frac{1+\beta}{2-\beta}}\right)
	\\&+
	2^{\frac{7\beta+7}{2-\beta}}\kappa^{\frac{2}{2-\beta}}(2\beta)C^{\frac{2}{2-\beta}}_{6}C^{\frac{4\beta}{2-\beta}}_{6\beta}
	\frac{b^{\frac{2}{2-\beta}}}{\nu^{\frac{\beta}{2-\beta}}}
	H^{\frac{\beta-1}{2-\beta}}\left(M_{2}^{\frac{1+\beta}{2-\beta}} + (4H)^{\frac{1+\beta}{2-\beta}}\right)
	+ \frac{\nu}{2l^{2}},\nonumber
\end{align}
Therefore
\begin{align}
	\frac{1}{2}\frac{d}{dt} \left(\|\nabla g(t)\|_{L^{2}}^{2}+ \frac{1}{l^{2}}\|g(t)\|_{L^{2}}^{2}\right) &\leq \left(\|\nabla g(t)\|_{L^{2}}^{2}+ \frac{1}{l^{2}}\|g(t)\|_{L^{2}}^{2}\right)\left(-\frac{3\eta}{8} + \frac{8\eta^{2}c_{0}h^{2}}{\nu} +\tilde{Z}_{1}+Z_{2}(t)\right)\nonumber\\&+ |\alpha-\beta|^{2} \left(Z_{3}+Z_{4}(t)\right)+ |\tilde{a}-\tilde{b}|^{2}\left(Z_{5}+Z_{6}(t)\right), \label{57fghfg}
\end{align}
for all $0\leq t <T^{*}$. Thus, using Gronwall's inequality,
\begin{align}
	\|\nabla g(t)\|_{L^{2}}^{2}&+ \frac{1}{l^{2}}\|g(t)\|_{L^{2}}^{2} \leq 
	e^{\left(-\frac{3\eta}{8} + \frac{8\eta^{2}c_{0}h^{2}}{\nu}+\tilde{Z}_{1}\right)t +\int_{0}^{t}Z_{2}(s)ds}\left(\|\nabla g(0)\|_{L^{2}}^{2}+ \frac{1}{l^{2}}\|g(0)\|_{L^{2}}^{2}\right)\nonumber\\&+
	2|\alpha-\beta|^{2} \int_{0}^{t}e^{\left(-\frac{3\eta}{8} + \frac{8\eta^{2}c_{0}h^{2}}{\nu}+\tilde{Z}_{1}\right)(t-r) +\int_{r}^{t}Z_{2}(s)ds}\left(Z_{3}+Z_{4}(r)\right)dr
	\nonumber\\&+
	2|\tilde{a}-\tilde{b}|^{2} \int_{0}^{t}e^{\left(-\frac{3\eta}{8} + \frac{8\eta^{2}c_{0}h^{2}}{\nu}+\tilde{Z}_{1}\right)(t-r) +\int_{r}^{t}Z_{2}(s)ds}\left(Z_{5}+Z_{6}(r)\right)dr, \nonumber
\end{align}
 $ 0\leq t < T^{*}$. Using again estimates given in Corollary \ref{da0101.4565.1}, we have
\begin{align}
	\left(-\frac{3\eta}{8} + \frac{8\eta^{2}c_{0}h^{2}}{\nu}+\tilde{Z}_{1}\right)(t-r) +\int_{r}^{t}Z_{2}(s)ds &\leq 
	\left(-\frac{3\eta}{8} + \frac{8\eta^{2}c_{0}h^{2}}{\nu}+\tilde{Z}_{1} + 2^{12}\kappa^{2}(2\beta)C^{2}_{6}C^{4\beta}_{6\beta}\frac{l}{\nu}M_{3}\right)(t-r) \nonumber\\&+ 2^{13}\kappa^{2}(2\beta)C^{2}_{6}C^{4\beta}_{6\beta}\frac{l}{\nu^{2}}M_{2}. \label{fbdf3ko}
\end{align}
Since, by hypothesis,
\begin{align}
	-\frac{3\eta}{8} + \frac{8\eta^{2}c_{0}h^{2}}{\nu}+\tilde{Z}_{1} + 2^{12}\kappa^{2}(2\beta)C^{2}_{6}C^{4\beta}_{6\beta}\frac{l}{\nu}M_{3}< -\frac{\eta}{8}, \label{fbdf3ko2}
\end{align}
it follows that
\begin{align}
	\|\nabla g(t)\|_{L^{2}}^{2}+ \frac{1}{l^{2}}\|g(t)\|_{L^{2}}^{2} &\leq 
	e^{-\frac{\eta}{8}t+ 2^{13}\kappa^{2}(2\beta)C^{2}_{6}C^{4\beta}_{6\beta}\frac{l}{\nu^{2}}M_{2}}\left(\|\nabla g(0)\|_{L^{2}}^{2}+ \frac{1}{l^{2}}\|g(0)\|_{L^{2}}^{2}\right)\nonumber\\&+
	2|\alpha-\beta|^{2} \int_{0}^{t}e^{-\frac{\eta}{8}(t-r) + 2^{13}\kappa^{2}(2\beta)C^{2}_{6}C^{4\beta}_{6\beta}\frac{l}{\nu^{2}}M_{2}}\left(Z_{3}+Z_{4}(r)\right)dr
	\nonumber\\&+
	2|\tilde{a}-\tilde{b}|^{2} \int_{0}^{t}e^{-\frac{\eta}{8}(t-r) + 2^{13}\kappa^{2}(2\beta)C^{2}_{6}C^{4\beta}_{6\beta}\frac{l}{\nu^{2}}M_{2}}\left(Z_{5}+Z_{6}(r)\right)dr
	\nonumber\\
	&\leq 
	e^{-\frac{\eta}{8}t+2^{13}\kappa^{2}(2\beta)C^{2}_{6}C^{4\beta}_{6\beta}\frac{l}{\nu^{2}}M_{2}}\left(\|\nabla g(0)\|_{L^{2}}^{2}+ \frac{1}{l^{2}}\|g(0)\|_{L^{2}}^{2}\right)\nonumber\\&+
	|\alpha-\beta|^{2} \frac{2^{4}}{\eta}e^{2^{13}\kappa^{2}(2\beta)C^{2}_{6}C^{4\beta}_{6\beta}\frac{l}{\nu^{2}}M_{2}}Z_{3}
	\nonumber\\&+
	2|\alpha-\beta|^{2} e^{2^{13}\kappa^{2}(2\beta)C^{2}_{6}C^{4\beta}_{6\beta}\frac{l}{\nu^{2}}M_{2}}\int_{0}^{t}e^{-\frac{\eta}{8}(t-r)}Z_{4}(r)dr
	\nonumber\\&+
	|\tilde{a}-\tilde{b}|^{2} \frac{2^{4}}{\eta}e^{2^{13}\kappa^{2}(2\beta)C^{2}_{6}C^{4\beta}_{6\beta}\frac{l}{\nu^{2}}M_{2}}Z_{5}
	\nonumber\\&+
	2|\tilde{a}-\tilde{b}|^{2} e^{2^{13}\kappa^{2}(2\beta)C^{2}_{6}C^{4\beta}_{6\beta}\frac{l}{\nu^{2}}M_{2}}\int_{0}^{t}e^{-\frac{\eta}{8}(t-r)}Z_{6}(r)dr,\nonumber
\end{align}
for all $0\leq t < T^{*}$. By Corollary \ref{da0101.4565.1.2}, we get
\begin{align}
	\int_{0}^{t}e^{-\frac{\eta}{8}(t-r)}Z_{4}(r)dr&= \frac{2^{14}\tilde{a}^{2}C^{10}_{10}}{(2-\max\{\alpha,\,\beta\})^{2}} \frac{l^{4}}{\nu^{7}}
	M_{2}^{4} \int_{0}^{t}e^{-\frac{\eta}{8}(t-r)}\|Au(r)\|^{2}_{L^{2}}dr \nonumber\\&\leq \frac{2^{14}\tilde{a}^{2}C^{10}_{10}}{(2-\max\{\alpha,\,\beta\})^{2}} \frac{l^{4}}{\nu^{7}}
	M_{2}^{4} \left(\frac{4}{\nu}M_{2}+ \frac{16}{\eta \nu^{\frac{3\alpha-2}{\alpha-1}} a^{\frac{1}{\alpha-1}}} M_{2} +\frac{32}{\eta \nu^{2}}\|f\|^{2}_{L^{\infty}_{t}L^{2}}\right), \label{577ggbdf}
\end{align}
and
\begin{align}
	\int_{0}^{t}e^{-\frac{\eta}{8}(t-r)}Z_{6}(r)dr&=  2^{11}C_{10}^{10}\frac{l^{4}}{\nu^{7}}M_{2}^{4}\int_{0}^{t}e^{-\frac{\eta}{8}(t-r)}\|Au(r)\|^{2}_{L^{2}}dr \nonumber\\&\leq 2^{11}C_{10}^{10}\frac{l^{4}}{\nu^{7}}M_{2}^{4} \left(\frac{4}{\nu}M_{2}+ \frac{16}{\eta \nu^{\frac{3\alpha-2}{\alpha-1}} a^{\frac{1}{\alpha-1}}} M_{2} +\frac{32}{\eta \nu^{2}}\|f\|^{2}_{L^{\infty}_{t}L^{2}}\right). \label{577ggbdf2}
\end{align}
Hence
\begin{align}
	&\|\nabla g(t)\|_{L^{2}}^{2}+ \frac{1}{l^{2}}\|g(t)\|_{L^{2}}^{2} 
	\leq 
	e^{2^{13}\kappa^{2}(2\beta)C^{2}_{6}C^{4\beta}_{6\beta}\frac{l}{\nu^{2}}M_{2}}\left\{
	e^{-\frac{\eta}{8}t}\left(\|\nabla g(0)\|_{L^{2}}^{2}+ \frac{1}{l^{2}}\|g(0)\|_{L^{2}}^{2}\right)\right.\nonumber\\&+
	2|\alpha-\beta|^{2} \left[\frac{2^{3}}{\eta}Z_{3}+
	\frac{2^{14}\tilde{a}^{2}C^{10}_{10}}{(2-\max\{\alpha,\,\beta\})^{2}} \frac{l^{4}}{\nu^{7}}
	M_{2}^{4} \left(\frac{4}{\nu}M_{2}+ \frac{16}{\eta \nu^{\frac{3\alpha-2}{\alpha-1}} a^{\frac{1}{\alpha-1}}} M_{2} +\frac{32}{\eta \nu^{2}}\|f\|^{2}_{L^{\infty}_{t}L^{2}}\right)\right]
	\nonumber\\&+\left.
	2|\tilde{a}-\tilde{b}|^{2} \left[\frac{2^{3}}{\eta}Z_{5}
	+2^{11}C_{10}^{10}\frac{l^{4}}{\nu^{7}}M_{2}^{4} \left(\frac{4}{\nu}M_{2}+ \frac{16}{\eta \nu^{\frac{3\alpha-2}{\alpha-1}} a^{\frac{1}{\alpha-1}}} M_{2} +\frac{32}{\eta \nu^{2}}\|f\|^{2}_{L^{\infty}_{t}L^{2}}\right)\right]\right\}<H \label{ftgh66556k}
\end{align}
for all $ 0\leq t < T^{*}$. Using continuity of $\|\nabla g(t)\|_{L^{2}}^{2}+ \frac{1}{l^{2}}\|g(t)\|_{L^{2}}^{2}$ and definition of $T^{*}$, we conclude that $T^{*}=\infty$.
Let $B,C$ and $D$ be the following quantities
\begin{align}
	B&= e^{2^{13}\kappa^{2}(2\beta)C^{2}_{6}C^{4\beta}_{6\beta}\frac{l}{\nu^{2}}M_{2}}, \label{zzzer1}\\\nonumber&\\
	C&= 2B\left\{\frac{2^{3}}{\eta}\left[ \frac{64\tilde{a}^{2}}{e^{2}}\frac{\nu}{l^{4}}M_{1} +
	\frac{2^{15}\tilde{a}^{2}(C^{12}_{6}+C^{10}_{10})}{(2-\max\{\alpha,\,\beta\})^{2}}\frac{l^{2}}{\nu^{7}}\left(M_{2}+\frac{1}{l^{2}}M_{1}\right)^{5}\right] \right.\nonumber\\& \left.+ 
	\frac{2^{14}\tilde{a}^{2}C^{10}_{10}}{(2-\max\{\alpha,\,\beta\})^{2}} \frac{l^{4}}{\nu^{7}}
	M_{2}^{4} \left(\frac{4}{\nu}M_{2}+ \frac{16}{\eta \nu^{\frac{3\alpha-2}{\alpha-1}} a^{\frac{1}{\alpha-1}}} M_{2} +\frac{32}{\eta \nu^{2}}\|f\|^{2}_{L^{\infty}_{t}L^{2}}\right)\right\}, \label{zzzer2}\\ \nonumber&\\
	D&= 2B\left\{\frac{2^{3}}{\eta}\left[6\frac{\nu}{l^{4}} M_{1}+
	2^{12}(C^{10}_{6}+ C_{10}^{10})\frac{l^{2}}{\nu^{7}}\left(M_{2}+\frac{1}{l^{2}}M_{1}\right)^{5}\right]
	\right.\nonumber\\&\left. +2^{11}C_{10}^{10}\frac{l^{4}}{\nu^{7}}M_{2}^{4} \left(\frac{4}{\nu}M_{2}+ \frac{16}{\eta \nu^{\frac{3\alpha-2}{\alpha-1}} a^{\frac{1}{\alpha-1}}} M_{2} +\frac{32}{\eta \nu^{2}}\|f\|^{2}_{L^{\infty}_{t}L^{2}}\right)\right\}. \label{zzzer3}
\end{align}
Therefore, the result asserted in Theorem \ref{CDA3a1.b.c.d} follows from (\ref{ftgh66556k}).

\begin{flushright}
	\rule{2mm}{2mm}
\end{flushright}


\section{Proof of Theorem \ref{CDA3a1.b.c.d.e}}\label{sec7} 


In the proof of Theorem \ref{CDA3a1.b.c.d}, we obtained that $T^{*}=\infty$.
Then, using Gronwall in (\ref{57fghfg}), we obtain
\begin{align}
	\|\nabla g(t)\|_{L^{2}}^{2}&+ \frac{1}{l^{2}}\|g(t)\|_{L^{2}}^{2} \leq 
	e^{\left(-\frac{3\eta}{8} + \frac{8\eta^{2}c_{0}h^{2}}{\nu}+\tilde{Z}_{1}\right)\left(t-\frac{2l^{2}}{\nu}\right) +\int_{\frac{2l^{2}}{\nu}}^{t}Z_{2}(s)ds}\left(\left\|\nabla g\left(\frac{2l^{2}}{\nu}\right)\right\|_{L^{2}}^{2}+ \frac{1}{l^{2}}\left\|g\left(\frac{2l^{2}}{\nu}\right)\right\|_{L^{2}}^{2}\right)\nonumber\\&+
	2|\alpha-\beta|^{2} \int_{\frac{2l^{2}}{\nu}}^{t}e^{\left(-\frac{3\eta}{8} + \frac{8\eta^{2}c_{0}h^{2}}{\nu}+\tilde{Z}_{1}\right)(t-r) +\int_{r}^{t}Z_{2}(s)ds}\left(Z_{3}+Z_{4}(r)\right)dr
	\nonumber\\&+
	2|\tilde{a}-\tilde{b}|^{2} \int_{\frac{2l^{2}}{\nu}}^{t}e^{\left(-\frac{3\eta}{8} + \frac{8\eta^{2}c_{0}h^{2}}{\nu}+\tilde{Z}_{1}\right)(t-r) +\int_{r}^{t}Z_{2}(s)ds}\left(Z_{5}+Z_{6}(r)\right)dr,\nonumber
\end{align}
for all $t \geq \frac{2l^{2}}{\nu}$, where $Z_{2}$, $Z_{3}$, $Z_{4}$, $Z_{5}$, $Z_{6}$ and $\tilde{Z}_{1}$ are given in (\ref{z2ul})-(\ref{z6ul}) and (\ref{dfbdfr54}).
Then, using (\ref{fbdf3ko}) and (\ref{fbdf3ko2}), we have
\begin{align}
	\|\nabla g(t)\|_{L^{2}}^{2}+ \frac{1}{l^{2}}\|g(t)\|_{L^{2}}^{2} &\leq
	e^{-\frac{\eta}{8}\left(t-\frac{2l^{2}}{\nu}\right)+2^{13}\kappa^{2}(2\beta)C^{2}_{6}C^{4\beta}_{6\beta}\frac{l}{\nu^{2}}M_{2}}\left(\left\|\nabla g\left(\frac{2l^{2}}{\nu}\right)\right\|_{L^{2}}^{2}+ \frac{1}{l^{2}}\left\|g\left(\frac{2l^{2}}{\nu}\right)\right\|_{L^{2}}^{2}\right)\nonumber\\&+
	|\alpha-\beta|^{2} \frac{2^{4}}{\eta}e^{2^{13}\kappa^{2}(2\beta)C^{2}_{6}C^{4\beta}_{6\beta}\frac{l}{\nu^{2}}M_{2}}Z_{3}
	\nonumber\\&+
	2|\alpha-\beta|^{2} e^{2^{13}\kappa^{2}(2\beta)C^{2}_{6}C^{4\beta}_{6\beta}\frac{l}{\nu^{2}}M_{2}}\int_{\frac{2l^{2}}{\nu}}^{t}e^{-\frac{\eta}{8}(t-r)}Z_{4}(r)dr
	\nonumber\\&+
	|\tilde{a}-\tilde{b}|^{2} \frac{2^{4}}{\eta}e^{2^{13}\kappa^{2}(2\beta)C^{2}_{6}C^{4\beta}_{6\beta}\frac{l}{\nu^{2}}M_{2}}Z_{5}
	\nonumber\\&+
	2|\tilde{a}-\tilde{b}|^{2} e^{2^{13}\kappa^{2}(2\beta)C^{2}_{6}C^{4\beta}_{6\beta}\frac{l}{\nu^{2}}M_{2}}\int_{\frac{2l^{2}}{\nu}}^{t}e^{-\frac{\eta}{8}(t-r)}Z_{6}(r)dr,\,\, \forall\, \nonumber
\end{align}
for all $t \geq \frac{2l^{2}}{\nu}$. By Lemma \ref{da0101.AA}, we get
\begin{align}
	\int_{\frac{2l^{2}}{\nu}}^{t}e^{-\frac{\eta}{8}(t-r)}Z_{4}(r)dr&= \frac{2^{14}\tilde{a}^{2}C^{10}_{10}}{(2-\max\{\alpha,\,\beta\})^{2}} \frac{l^{4}}{\nu^{7}}
	M_{2}^{4} \int_{\frac{2l^{2}}{\nu}}^{t}e^{-\frac{\eta}{8}(t-r)}\|Au(r)\|^{2}_{L^{2}}dr \leq \frac{2^{14}\tilde{a}^{2}C^{10}_{10}}{(2-\max\{\alpha,\,\beta\})^{2}} \frac{l^{4}}{\nu^{7}}
	M_{2}^{4} \frac{8}{\eta}M_{8}^{2}; \nonumber \\&\nonumber\\
	\int_{\frac{2l^{2}}{\nu}}^{t}e^{-\frac{\eta}{8}(t-r)}Z_{6}(r)dr&=  2^{11}C_{10}^{10}\frac{l^{4}}{\nu^{7}}M_{2}^{4}\int_{\frac{2l^{2}}{\nu}}^{t}e^{-\frac{\eta}{8}(t-r)}\|Au(r)\|^{2}_{L^{2}}dr \leq 2^{11}C_{10}^{10}\frac{l^{4}}{\nu^{7}}M_{2}^{4} \frac{1}{\eta}M_{8}^{2}. \nonumber
\end{align}
Consequently,
\begin{align}
	\|\nabla g(t)\|_{L^{2}}^{2}+ \frac{1}{l^{2}}\|g(t)\|_{L^{2}}^{2} 
	&\leq 
	e^{2^{13}\kappa^{2}(2\beta)C^{2}_{6}C^{4\beta}_{6\beta}\frac{l}{\nu^{2}}M_{2}}\left\{
	e^{-\frac{\eta}{8}\left(t-\frac{2l^{2}}{\nu}\right)}\left(\left\|\nabla g\left(\frac{2l^{2}}{\nu}\right)\right\|_{L^{2}}^{2}+ \frac{1}{l^{2}}\left\|g\left(\frac{2l^{2}}{\nu}\right)\right\|_{L^{2}}^{2}\right)\right.\nonumber\\&+
	2|\alpha-\beta|^{2} \left[\frac{2^{3}}{\eta}Z_{3}+
	\frac{2^{14}\tilde{a}^{2}C^{10}_{10}}{(2-\max\{\alpha,\,\beta\})^{2}} \frac{l^{4}}{\nu^{7}}
	M_{2}^{4} \frac{1}{\eta}M_{8}^{2}\right]
	\nonumber\\&+\left.
	2|\tilde{a}-\tilde{b}|^{2} \left[\frac{2^{3}}{\eta}Z_{5}
	+2^{11}C_{10}^{10}\frac{l^{4}}{\nu^{7}}M_{2}^{4}\frac{1}{\eta}M_{8}^{2} \right]\right\}, \label{ftgh66556k.3}
\end{align}
for all $t \geq \frac{2l^{2}}{\nu}$. Let $B$ be given in (\ref{zzzer1}), and
\begin{align}
	\tilde{C}&= \frac{2B}{\eta}\left\{2^{3}\left[ \frac{64\tilde{a}^{2}}{e^{2}}\frac{\nu}{l^{4}}M_{1} +
	\frac{2^{15}\tilde{a}^{2}(C^{12}_{6}+C^{10}_{10})}{(2-\max\{\alpha,\,\beta\})^{2}}\frac{l^{2}}{\nu^{7}}\left(M_{2}+\frac{1}{l^{2}}M_{1}\right)^{5}\right] \right.\nonumber\\& \left.+ 
	\frac{2^{14}\tilde{a}^{2}C^{10}_{10}}{(2-\max\{\alpha,\,\beta\})^{2}} \frac{l^{4}}{\nu^{7}}
	M_{2}^{4} M^{2}_{8}\right\}, \label{zzzer2.2}\\ \nonumber&\\
	\tilde{D}&= \frac{2B}{\eta}\left\{2^{3}\left[6\frac{\nu}{l^{4}} M_{1}+
	2^{12}(C^{10}_{6}+ C_{10}^{10})\frac{l^{2}}{\nu^{7}}\left(M_{2}+\frac{1}{l^{2}}M_{1}\right)^{5}\right]
	\right.\nonumber\\&\left. +2^{11}C_{10}^{10}\frac{l^{4}}{\nu^{7}}M_{2}^{4}M_{8}^{2}\right\}. \label{zzzer3.3}
\end{align}
Finally, the result stated in Theorem \ref{CDA3a1.b.c.d.e} follows from (\ref{ftgh66556k.3}).

\begin{flushright}
	\rule{2mm}{2mm}
\end{flushright}

\section{Conclusion}
In this paper, we proved that under suitable conditions, it is possible to approximate physical solutions of the three-dimensional Brinkman-Forchheimer-extended Darcy model even with the parameters related to the damping term unknown. Computational experiments, as well as analysis of algorithms to recovery the physical parameters $\alpha$ and $a$ based in continuous data assimilation techniques will be subject of forthcoming work.




\begin{thebibliography}{99}
	
	
	\bibitem{Albanez1} \textsc{Albanez, D.A.F., Nussenzveig Lopes, H., Titi, E.} \textit{Continuous data assimilation for the three-dimensionl Navier-Stokes-$\alpha$ model}, Asymptotic Analysis, 97, 139-164, 2016.

    \bibitem{Albanez2} \textsc{Albanez, D.A.F., Benvenutti, M.J.} \textit{Continuous data assimilation algorithm for simplified Bardina model.}, Evolution Equations \& Control Theory, 7.1, 33-52, 2018.
	

	
	\bibitem{Amao} \textsc{Amao A. M.} \textit{ Mathematical model for Darcy Forchheimer flow with applications to well performance analysis},  Thesis in petroleum engineering for the Degree of MASTER OF SCIENCE, Faculty of Texas Tech University, 2007.
	
		\bibitem{Azouani} \textsc{Azouani A., Olson E., Titi E.} \textit{Continuous data assimilation using General interpolant observables}, J Nonlinear Sci., 24 , 277-304, 2014.
	
	\bibitem{Azouani2} \textsc{Azouani A., Titi E.}  \textit{Feedback control of nonlinear dissipative systems by finite determining parameters - a reaction-diffusion paradigm}, Evolution Equations \& Control Theory, 3(4), (2014).
	
	
	\bibitem{Biswas} \textsc{Biswas A.,  Hudson J.,  Larios A., Pei Y.} \textit{Continuous data assimilation for the 2D magnetohydrodynamic equations using one component of the velocity and magnetic field}, Asymptotic Analysis 0, 1–43, 2017.
	
	\bibitem{Biswas2} \textsc{Biswas A.,  Hudson J.,  Larios A., Pei Y.} \textit{Continuous Data Assimilation for the Three Dimensional Navier-Stokes Equations}, SIAM J. Math. Anal., 53(6), 6697–6723, 2021.
	
	\bibitem{Cai} \textsc{Cai X., Jiu Q.} \textit{Weak and strong solutions for the incompressible Navier-Stokes equations with damping}, J. Math Analysis and applications, 343 , 799-809, 2008.
	
	\bibitem{Carlson1} \textsc{Carlson E., Hudson J., Larios A.} \textit{Parameter Recovery for the 2 Dimensional Navier-Stokes Equations via Continuous Data Assimilation}, SIAM Journal on Scientific Computing, 42(1), A250-A270, 2020.
	
	\bibitem{Carlson2} \textsc{Carlson E., Hudson J., Larios A., Martinez V.R., Ng E., Whitehead J.} \textit{Dynamically learning the parameters of a chaotic system using partial observations}, Discrete and Continuous Dynamical Systems, 48(8), 3809-3839, 2021.
	
	\bibitem{Farhat1} \textsc{Farhat A., Jolly M., Titi E.} \textit{Continuous data assimilation for the 2D B\'enard convection through velocity measurements}, Physica D,. 303, 59-66, 2015.
	
	\bibitem{Farhat2} \textsc{Farhat A., Lunasin E., Titi E.} \textit{Abridged continuous data assimilation for the 2D Navier-Stokes equations utilizing measurements of only one component of the velocity field}, J.Math. Fluid Mechanics, 18, 1-23, 2016.
	
	\bibitem{Farhat4} \textsc{Farhat A., Lunasin E., Titi E.} \textit{Data assimilation algorithm for 3D B\'enard convection in porous media emplying only temperature measurements}, J. Math. Analysis and Applications, 438, 492-506, 2016.
	
	\bibitem{Rosa} \textsc{Foias C., Manley O., R. Rosa, Temam R.} \textit{%
		Navier-Stokes equations and turbulence}, Cambridge, 2004.
	
	\bibitem{Friedman} \textsc{Friedman A.} \textit{Partial differential equations}, Dover Publications, 2008.
	
    \bibitem{Jolly1} \textsc{Jolly, M.S., Martinez, V.R., Titi, E.S.} \textit{A data assimilation algorithm for the subcritical surface
	quasi-geostrophic equation.}, Advanced Nonlinear Studies, 17.1, 167–192, 2017.

    \bibitem{Jolly2} \textsc{Jolly, M.S., Sadigov, T., Titi, E.S.} \textit{A determining form for the damped driven nonlinear Schrödinger
	equation—Fourier modes case.}, Journal of Differential Equations, 258.8, 2711–2744, 2015.
 
    \bibitem{Jolly3} \textsc{Jolly, M.S., Sadigov, T., Titi, E.S.} \textit{Determining form and data assimilation algorithm for weakly damped and driven Korteweg–de Vries equation—Fourier modes case.}, Nonlinear Analysis: Real World Applications, 36, 287-317, 2017.
 
	\bibitem{Giorgi} \textsc{Giorgi T.} \textit{Derivation of the Forchheimer law via matched
		asymptotic expansions}, Transport in Porous Media 29: 191-206, 1997.
	
	\bibitem{HSU} \textsc{Hsu C. T., Cheng P,} \textit{Thermal dispersion in a porous medium}, Int. J. Heat Mass Transfer Vol. 33. No. 8., 1587-177. 1990.
	
	\bibitem{Ingham} \textsc{Ingham D. B., Pop I.} \textit{Transport Phenomena in Porous Media II}, Pergamon 2002.
	
	\bibitem{Joseph} \textsc{Joseph D. D., Nield D. A., Papanicolaou G.} \textit{Nonlinear Equation Governing Flow in a Saturated Porous Medium}, Water Resources Research,  18,  1049-1052, 1982.
	
	
	\bibitem{Li} \textsc{Kim Y., Li K.} \textit{Time-periodic strong solutions of the 3D Navier-Stokes equations with damping}, Electronic Journal of  Differential equations, 244, 1-11, 2017.
	
	\bibitem{Li2} \textsc{Kim Y., Li K., Kim C.} \textit{Uniqueness and regularity for the 3D Boussinesq system with damping}, Annali Dell'Univertita'Di Ferrara, 67, 149-173, 2021.
	
	\bibitem{Nield} \textsc{Lage, J. L., Antohe B. V., Nield. D. A. } \textit{Two types of nonlinear
		pressure-drop versus flow-rate relation observed for saturated porous media}, Journal of Fluid Mechanics J. Fluids Eng, 700-706, 1997.
	
	\bibitem{Markowich} \textsc{Markowich P., Titi E. \& Trabelsi S.} \textit{Continuous data assimilation for the three-dimensional Brinkman–Forchheimer-extended Darcy model}, Nonlinearity, Volume 29(4), 1291-1328, 2016.
	
	\bibitem{Mei} \textsc{Mei C., Auriault J.L. } \textit{Upscaling Forchheimer law}, Transport in Porous Media, 70, 2013-229, 2007.
	
	\bibitem{Mei2} \textsc{Mei C., Auriault J.L. } \textit{The effect of weak inertia on flow through a porous medium}, Journal of Fluid Mechanics 222, 647-663, 1991.
	
	\bibitem{Nield2} \textsc{Nield D. A., Bejan A.} \textit{Convection in Porous Media, Third edition}, Springer, 2006.
	
	\bibitem{Pardo} \textsc{Pardo D., Valero J. \& Gim\'enez A.} \textit{Global attractors for weak solutions of the three-dimensional Navier-Stokes equations with damping}, Discrete e Continuous Dynamical Systems  24 (8), 3569-3590, 2019.
	
	\bibitem{Skjetne} \textsc{Skjetne E., Auriault J.L.} \textit{High-Velocity Laminar and Turbulent
		Flow in Porous Media}, Transport in Porous Media 36: 131–147, 1999.
	
	\bibitem{Temam} \textsc{Temam R.} \textit{Navier-Stokes Equations: Theory
		and Numerical Analysis}, AMS Chelsea Publishing, 2001.
	
	\bibitem{Temam2} \textsc{Temam R.} \textit{Navier-Stokes Equations and Nonlinear Functional Analysis}, Society for Ind. Appl. Math., 1995.
	
	\bibitem{Vafai} \textsc{Vafai K., Kim, S.J.} \textit{On the limitations of the Brinkman-Forchheimer-extended Darcy equation}, Int. J.  Heat and Fluid Flow, 16,1-15, 1995.
	
	\bibitem{Zhou} \textsc{Wang W., Zhou G.} \textit{Remarks on the regularity criterion of the Navier-Stokes equations with nonlinear damping}, Mathematical Problems in Engineering, 1-5, 2015.
	
	\bibitem{Whitaker} \textsc{Whitaker S.} \textit{The Forchheimer Equation: A Theoretical
		Development}, Transport in Porous Media 25: 27-61, 1996.
	
	\bibitem{Lu} \textsc{Zhang Z., Wu X., Lu M.} \textit{On the uniqueness of strong solution to the incompressible Navier-Stokes equations with damping}, J. Math Analysis and applications, 377, 414-419, 2011.
	
	\bibitem{Zhong} \textsc{Zhong X.} \textit{A note on the uniqueness of strong solution to the incompressible Navier-Stokes equations with damping}, Electronic Journal of qualitative Theory of Differential equations, 15, 1-4, 2019.
	
	\bibitem{Zhong2} \textsc{Zhou Y.} \textit{Regularity and uniqueness for the 3D incompressible Navier-Stokes equations with damping}, Applied mathematics Letters, 25, 1822-1825, 2012.
	

	
\end{thebibliography}
\end{document}